\documentclass[12pt]{scrartcl}

\usepackage{amsmath,amssymb,amsthm,mathtools,array}
\usepackage{graphicx}

\mathtoolsset{showonlyrefs=true}

\newcounter{thm}%[section]
\newenvironment{thm}[1][]{\refstepcounter{thm}\par\medskip
   \noindent \textbf{\sffamily Theorem~\thethm. #1} \slshape }{\medskip}

\newcounter{rem}%[section]
\newenvironment{rem}[1][]{\refstepcounter{rem}\par\medskip
   \noindent \textbf{\sffamily Remark~\therem. #1} \rmfamily }{\medskip}

\usepackage{tikz}

\usetikzlibrary{patterns}

\newlength{\hatchspread}
\newlength{\hatchthickness}
\newlength{\hatchshift}
\newcommand{\hatchcolor}{}
\tikzset{hatchspread/.code={\setlength{\hatchspread}{#1}},
         hatchthickness/.code={\setlength{\hatchthickness}{#1}},
         hatchshift/.code={\setlength{\hatchshift}{#1}},% must be >= 0
         hatchcolor/.code={\renewcommand{\hatchcolor}{#1}}}
\tikzset{hatchspread=3pt,
         hatchthickness=0.4pt,
         hatchshift=0pt,% must be >= 0
         hatchcolor=black}
\pgfdeclarepatternformonly[\hatchspread,\hatchthickness,\hatchshift,\hatchcolor]% variables
   {custom north west lines}% name
   {\pgfqpoint{\dimexpr-2\hatchthickness}{\dimexpr-2\hatchthickness}}% lower left corner
   {\pgfqpoint{\dimexpr\hatchspread+2\hatchthickness}{\dimexpr\hatchspread+2\hatchthickness}}% upper right corner
   {\pgfqpoint{\dimexpr\hatchspread}{\dimexpr\hatchspread}}% tile size
   {% shape description
    \pgfsetlinewidth{\hatchthickness}
    \pgfpathmoveto{\pgfqpoint{0pt}{\dimexpr\hatchspread+\hatchshift}}
    \pgfpathlineto{\pgfqpoint{\dimexpr\hatchspread+0.15pt+\hatchshift}{-0.15pt}}
    \ifdim \hatchshift > 0pt
      \pgfpathmoveto{\pgfqpoint{0pt}{\hatchshift}}
      \pgfpathlineto{\pgfqpoint{\dimexpr0.15pt+\hatchshift}{-0.15pt}}
    \fi
    \pgfsetstrokecolor{\hatchcolor}
    \pgfusepath{stroke}
   }

\newcommand\gp{\,}

\begin{document}

\begin{center}
\bfseries{\sffamily\Large
Paravectors and the Geometry of 3D Euclidean Space}
\end{center}

\bigskip
\begin{center}
{\bfseries\large Jayme Vaz Jr.}\footnote{ On leave of absence 
from Department of Applied Mathematics, IMECC, 
University of Campinas, 
Campinas, SP, Brazil}$^,$\footnote{\texttt{<vaz@uwaterloo.ca>, <vaz@ime.unicamp.br>}}
{\bfseries\large and} {\bfseries\large Stephen Mann}\footnote{\texttt{<smann@uwaterloo.ca>}}\\
\medbreak
\small 
Cheriton School of Computer Science\\
University of Waterloo\\
Waterloo, ON, Canada
\end{center}
\normalsize

\bigskip

\footnotesize 

\begin{center}
\begin{minipage}{11 cm}
\textbf{Abstract:} 
We introduce the concept of paravectors to describe the
geometry of points in a three dimensional space. After defining a 
suitable product of paravectors, we introduce the concepts of 
biparavectors and triparavectors to describe line
segments and plane fragments in this space. 
A key point in this product of paravectors 
is the notion of the orientation of a point, in such a way that 
biparavectors representing line segments are the result of
the product of points with opposite orientations. Incidence 
relations can also be formulated in terms of the product 
of paravectors. 
To study the transformations of points, lines, and planes, we 
introduce an algebra of transformations that is analogous to the 
algebra of creation and annihilation operators in quantum 
theory. The paravectors, biparavectors and triparavectors are
mapped into this algebra and their transformations are 
studied; we show that this formalism describes in an 
unified way the operations of reflection, rotations (circular 
and hyperbolic), translation, shear and non-uniform 
scale transformation. Using the concept of Hodge duality, 
we define a new operation called cotranslation, and show
that the operation of perspective projection can be 
written as a composition of the translation and cotranslation
operations.  We also show that the operation of pseudo-perspective can be
implemented using the cotranslation operation.

\medskip


\medskip

\end{minipage}
\end{center}

\normalsize

\section{Introduction}

Geometric transformations are part of the language
we use to ground our scientific knowledge, and consequently 
their study is of paramount importance. 
Two examples illustrate this point.  From an abstract 
point of view, ever since Klein's Erlangen 
Program~\cite{Klein1,Klein2}, geometric transformations
are considered as part of the geometry itself, so that 
a particular kind of geometry consists of a space of
objects and their transformation group; from a practical 
point of view, an entire area such as computer graphics 
has geometric transformations as the basis of the
modelling of its objects.
The objective of 
this work is to study the geometric transformations 
used in computer graphics and to develop an algebraic formalism 
for their description that reflects the geometry used in computer graphics.

The usual mathematical background of  
computer graphics is linear algebra. Objects are
represented by vectors and their motions by linear
transformations. Given the choice of a basis, vectors
are represented by column (or row) matrices and
linear transformations by matrices. Although from 
the computational point of view the use of matrices
is a natural option, from the theoretical point of
view the use of matrices is not a good choice because coordinates
do not have material existence. Ideally,
objects and their motions should be described through 
algebraic relations involving abstract vectors, leaving
the use of coordinates as a final step after choosing
an arbitrary basis and origin.
Several authors describe the use of quaternions 
in computer graphics~\cite{Conway,Shoemake,Goldman,Hanson,Kuipers,Vince3}. 
Although quaternions are elegant, their attractiveness diminishes
when one has to deal with translations of points, 
and at this point a quaternion 
enthusiast usually returns to
a matrix formulation (with dual quaternions being a less common
approach that handles translation~\cite{Kenwright}). 
However, quaternions are only one particular
Clifford algebra~\cite{VazRocha,DFM}, and
it makes sense to consider other Clifford algebras that 
incorporate transformations not possible with quaternions.

Two well-known alternative Clifford algebras
are the homogeneous model, which uses
vectors in a 4D vector space to 
represent points in a 3D space, and the
conformal model, which uses the conformal compactification  
of the 3D space in a 5D space to describe translations and conformal 
transformations~\cite{DFM}.
However, for computer graphics, neither model is
fully satisfactory because the  non-conformal transformations used
in the graphics pipeline, like shear or
non-uniform scaling, cannot be described in these models as rotors.

The primary realization is that computer graphics uses affine spaces~\cite{Goldman2000,Goldman2002} together with a single projection, and other Clifford algebra models have been proposed
for computer graphics, including~\cite{GS} and~\cite{Dorst2}. 
The model developed Goldman-Mann~\cite{GS} studies the standard
transformations of points used in computer graphics by means 
of the Clifford algebra of the vector space $\mathbb{R}^{4,4}$, 
while the model of Dorst~\cite{Dorst2} focuses on
transformations of lines using the Clifford algebra
of the vector space $\mathbb{R}^{3,3}$. Although 
the model of Goldman-Mann~\cite{GS} can describe
all the standard transformations in a similar way, 
their model is based on an 8 dimensional vector space.
The main motivation of our work is to study the 
standard transformations of computer graphics 
in a somewhat analogous way as done
by Goldman-Mann~\cite{GS} but with a smaller algebraic structure.

To accomplish our goal, we need to look 
for a suitable algebraic structure for describing
the objects of a 3D affine space 
and to describe the transformations of the objects 
in this space by means of endomorphisms of this algebraic
structure. Since Clifford algebras are subalgebras of
the algebra of endomorphisms of the exterior algebra, 
we will start looking to the exterior algebra of a 
3D vector space. 
The central object of our model is the paravector.
A paravector is an object composed of a 
scalar part and a vector part~\cite{Porteous}. 
It has been successfully used in alternative formulations 
of the special theory of relativity~\cite{Baylis1,Baylis2,Baylis3} 
and relativistic quantum mechanics~\cite{Baylis4,Vaz2016,Vaz2018} 
with a smaller algebraic structure than their 
traditional formulations. This fact led
us to the idea of describing 
the objects of the 3D affine space 
using a model based on paravectors.

Our plan is to describe points by means of 
paravectors, and the objects constructed 
from points by some product of paravectors. 
If we know how points can be described
by paravectors, then it is natural to think 
of line segments as given by a product 
of two paravectors, which results in a 
new object that we call a $2$-paravector. 
Continuing with this reasoning, it is 
natural to think of plane fragments as 
described by the product of three paravectors, 
defining a new object called a $3$-paravector. 
This product is similar to the
exterior product of vectors~\cite{VazRocha,ext_alg}. 
However, as we will see, the product of paravectors 
is not the usual exterior product of vectors, 
but a version such that, in analogy to the
usual exterior product, the product of
a $k$-paravector and an $l$-paravector gives a
$(k+l)$-paravector. Notwithstanding, this 
description of geometric 
objects does not exhaust the problem, for 
we must know how to describe their geometric 
transformations. The remarkable fact 
is that from operators constructed from 
the products of paravectors, we will be able to  
describe several geometric 
transformations through algebraic transformations 
of the form $x \mapsto U x \bar{U}$. 
All this modeling will be done based
on a 3D vector space. 

However, given that our starting point is 
Grassmann's exterior algebra, our approach 
is in no way restricted to three dimensions---there
is no difficultly to generalize our 
approach to $n$-dimensions.  Moreover, 
although we will work with an orthogonal basis, 
the $\mathbb{Z}_n$-graded structure of the 
exterior algebra does not depends on this 
fact, and our approach can be generalized 
to arbitrary basis; regardless, we work with 
an orthogonal basis to avoid unnecessary 
complications. At this point it is worth  
remembering that Clifford algebras are $\mathbb{Z}_2$-graded 
algebras, and they inherit the $\mathbb{Z}_n$-graded 
structure of exterior algebras only when one works 
with an orthogonal basis. 
Nevertheless, complications are expected 
in some calculations because some 
algebraic manipulations can be done more easily
in terms of Clifford algebra than in Grassmann's exterior algebra. 
Regardless, since we want to 
model things in terms of a basic and general structure, 
we think that the
structure of Grassmann's exterior algebra is the appropriate
structure for a first approach.

We have organized this work as follows. 
In Section~\ref{sec.2} we give a brief introduction 
to the exterior algebra of the vector space 
$\mathbb{R}^3$ and some structures that
are defined on this exterior algebra. The exterior algebra 
of $\mathbb{R}^3$ is of
fundamental importance for this work, 
so we need to dedicate a few words to it, 
and also to set our notation. 
We define the exterior product of vectors,  
the multivector structure of the algebra, 
some algebraic operations on the exterior
algebra, and when $\mathbb{R}^3$ is endowed
with a scalar product, we define the
interior product and the Hodge star operator 
on the exterior algebra.  

In Section~\ref{sec.3} we introduce the concept
of paravectors, and, in a more general way, 
of a $k$-paravector, and then we define a 
product of paravectors based on the exterior 
product. To use this product to define 
a product of paravectors with a geometric 
interpretation, we introduce the idea of 
orientation of a point. The paravector 
representation of a point reminds us 
of the geometry of mass points introduced
by M\"obius~\cite{DFM,Coxeter,Goldman3}, and which we prefer to call
weighted points. We interpret the absolute value of
the scalar part of a paravector as the weight of a point, 
while the sign of the paravector is interpreted as the orientation
of the point. With this interpretation, the product 
of paravectors with opposite orientations provides 
a representation of a line segment in terms 
of a $2$-paravector that resembles
the representation of lines in terms of Pl\"ucker
coordinates~\cite{DFM}. Orientation is a 
critical issue for geometric computations~\cite{Stolfi}, 
and we incorporate orientation in the basis of our approach. 

Given the geometric objects, the next step is to
study their transformations. To investigate transformations, we
first note that the space of paravectors, and $k$-paravectors in general, are
subspaces of the vector space underlying the exterior algebra. 
Therefore, operations on elements of the exterior algebra 
can also act on $k$-paravectors. However, if we consider
left and right exterior products and left and right 
interior products on $k$-paravectors, we do not have
an associative structure, which, from the point of
view of computations, is inconvenient.
But we can get around this situation. If we look to the
exterior and interior products as operators, it is
known that these operators satisfy an algebra
analogous to the algebra of the creation and annihilation 
operators of fermions in Quantum Field Theory~\cite{VazRocha,QFT}. 
We review this approach in Section~\ref{sec.4}. The idea then
is to map a $k$-paravector into the algebra of 
transformations and work with the operator image 
of the $k$-paravector instead of the $k$-paravector
itself. The advantage of this procedure is 
that we have an associative
linear structure similar to the algebra of matrices, 
and which is suitable for the description 
of geometric transformations; as usual, this approach involves
the use of a dual
space.

The algebraic structure discussed in Section~\ref{sec.4} 
is used to describe some transformations of points
in Section~\ref{sec.5}. We discuss the reflection 
of points in a plane; non-uniform scale 
and shear transformations of points; rotation 
and hyperbolic rotation of points; and translation of
points. We also define a new transformation that 
we call cotranslation, and show that a 
composition of translation and cotranslation of a point 
gives the perspective projection of this point 
from the eye into the  
perspective plane, and we show that cotranslation of a point can give
pseudo-perspective.
Then in Section~\ref{sec.6} 
we discuss how to extend these transformations to
the $2$-paravectors and $3$-paravectors 
to describe transformations of lines and planes. 
Finally, in Section~\ref{sec.7} we present our 
conclusions, where we also compare our model for affine geometry
and perspective projections to two other Clifford algebras that model
similar things.

\begin{rem}
We end this introduction with a remark about some notation 
used in this work because of the different meaning of some font 
typefaces. Points are denoted by $P$, $Q$, etc.,
their representation in terms of paravectors (Section~\ref{sec.3}) 
by $\mathsf{P}$, $\mathsf{Q}$, etc., and their operator representation  
(Section~\ref{sec.5}) by $\mathbf{P}$, $\mathbf{Q}$, etc.
Vectors are denoted by $\vec{p}$, $\vec{q}$, etc., and their
operator representation by $\mathbf{p}$, $\mathbf{q}$, etc. 
Sometimes we also use the notation $\overrightarrow{PQ} = \vec{q} - \vec{p}$. 
\end{rem}

\section{The Exterior Algebra}
\label{sec.2}

Let $\mathbb{R}^3$ be a three dimensional vector space over $\mathbb{R}$ with  
basis $\mathfrak{B} = \{\vec{e}_1,\vec{e}_2,\vec{e}_3\}$, and let 
$(\bigwedge(\mathbb{R}^3),\wedge)$ be its exterior algebra, where 
$\wedge$ denotes the exterior product $\vec{v}\wedge\vec{u} = - \vec{u}\wedge\vec{v}$, 
$\forall \vec{v},\vec{u} \in \mathbb{R}^3$,   and 
$\bigwedge(\mathbb{R}^3)$ is defined as 
\begin{equation}
\label{ext.space.V}
\textstyle{\bigwedge}(\mathbb{R}^3) = \textstyle{\bigwedge^0}(\mathbb{R}^3) \oplus \textstyle{\bigwedge^1}(\mathbb{R}^3) \oplus 
\textstyle{\bigwedge^2}(\mathbb{R}^3) \oplus \textstyle{\bigwedge^3}(\mathbb{R}^3),
\end{equation}
where $\textstyle{\bigwedge^k}(\mathbb{R}^3)$ denotes the
vector space of $k$-vectors ($k=0,1,2,3$) with 
$\bigwedge^0(\mathbb{R}^3) = \mathbb{R}$ 
and $\bigwedge^1(\mathbb{R}^3) = \mathbb{R}^3$. 
An arbitrary element $\Phi \in \bigwedge(\mathbb{R}^3)$ is called a multivector. 
For more details, see \cite{VazRocha}. 

There are three important operations, which are involutions (i.e., transformations whose square is the identity), 
that can be defined on $\bigwedge(\mathbb{R}^3)$. Given $A_k \in \bigwedge^k(\mathbb{R}^3)$  
{\em grade involution} (or parity), denoted by a hat, is defined as 
\begin{equation}
\hat{A}_k = (-1)^k A_k ,
\end{equation}
{\em reversion}, denoted by a tilde, is defined as 
\begin{equation}
\label{eq.tilde}
\tilde{A}_k = (-1)^{k(k-1)/2} A_k , 
\end{equation}
and {\em conjugation}, denoted by a bar, is the composition of 
reversion and grade involution, 
\begin{equation}
\bar{A}_k = \widehat{\tilde{A}}_k = \widetilde{\hat{A}}_k = 
(-1)^{k(k+1)/2} A_k .  
\end{equation}
Important properties of these operations are 
\begin{equation}
\widehat{(A_k \wedge B_l)} = \widehat{A}_k \wedge \widehat{B}_l , 
\quad 
\widetilde{(A_k \wedge B_l)} = \widetilde{B}_l \wedge \widetilde{A}_k , 
\quad 
\overline{(A_k \wedge B_l)} = \bar{B}_l \wedge \bar{A}_k .
\end{equation}

We will denote the projectors $\bigwedge(\mathbb{R}^3) \rightarrow 
\bigwedge^k(\mathbb{R}^3)$ by $\langle \quad \rangle_k$, that is, 
if $A_k \in \bigwedge^k(\mathbb{R}^3)$ and 
\begin{equation}
A = A_0 + A_1 + A_2 + A_3, 
\end{equation}
then 
\begin{equation}
\langle A \rangle_k = A_k . 
\end{equation}
The effect of the three involutions on $A$ is 
\begin{gather}
\hat{A} = A_0 - A_1 + A_2 - A_3 , \\
\tilde{A} = A_0 + A_1 - A_2 - A_3 , \\
\bar{A} = A_0 - A_1 - A_2 + A_3 . 
\end{gather}

\subsection{The Interior Product}

Suppose that $\mathbb{R}^3$ is endowed with a
scalar product $(\vec{v}|\vec{u})$. Then we write 
\begin{equation}
g_{ij} = (\vec{e}_i|\vec{e}_j) , \quad i,j = 1,2,3. 
\end{equation}
We do not need to suppose that the basis $\mathfrak{B}$ is an orthonormal 
basis; the formalism developed in this work is general and 
can work with any basis. However, in general the use of orthonormal
bases simplifies many calculations, and is therefore a useful assumption. 
To avoid unnecessary difficulties, let us assume therefore
that the basis $\{\vec{e}_1,\vec{e}_2,\vec{e}_3\}$ 
is an orthonormal basis: 
\begin{equation}
g_{ij} = \delta_{ij} . 
\end{equation}
This orthonormal basis does not need, of course, to be
associated with Cartesian coordinates, but can be
associated with any orthogonal curvilinear coordinates. 

Our objective now is to extend this product, from $\mathbb{R}^3$, to 
$\bigwedge(\mathbb{R}^3)$.  First, let us fix a vector $\vec{v}$. Then 
we can consider the operation $\vec{u} \mapsto 
(\vec{v}|\vec{u})$ as a product that takes 
$\vec{u} \in \mathbb{R}^3 = \bigwedge^1(\mathbb{R}^3)$ to 
$(\vec{v}|\vec{u}) \in \mathbb{R} = \bigwedge^0(\mathbb{R}^3)$, 
and consider its generalization to $\bigwedge(\mathbb{R}^3)$ as 
a product that takes an element of $\bigwedge^k(\mathbb{R}^3)$ and gives
an element of $\bigwedge^{k-1}(\mathbb{R}^3)$ for $k = 0,1,2,3$. 
At this point it is useful to use a different notation, that is, 
we will denote this product by a dot, so that
\begin{equation}
\vec{v} \cdot \vec{u} = (\vec{v}|\vec{u}) . 
\end{equation}
In the case of scalars, we define 
\begin{equation}
\label{vacuum}
\vec{v}\cdot 1 = 0 , 
\end{equation}
for bivectors, 
\begin{equation}
\label{int.prod.biv}
\vec{v}\cdot (\vec{u}\wedge \vec{w}) = 
(\vec{v}\cdot \vec{u})\vec{w} - 
(\vec{v}\cdot \vec{w})\vec{u} , 
\end{equation}
and for trivectors, 
\begin{equation}
\vec{v}\cdot(\vec{u}\wedge\vec{w}\wedge \vec{z}) = 
(\vec{v}\cdot\vec{u}) \vec{w}\wedge \vec{z} - 
(\vec{v}\cdot\vec{w}) \vec{u}\wedge\vec{z} + 
(\vec{v}\cdot\vec{z})\vec{u}\wedge\vec{w} .
\end{equation}
From the last two equations, note that 
\begin{equation}
\label{aux}
\begin{split}
\vec{v}\cdot(\vec{u}\wedge\vec{w}\wedge \vec{z}) & = 
(\vec{v}\cdot(\vec{u}\wedge\vec{w})) \wedge \vec{z} + 
\vec{u}\wedge\vec{w} (\vec{v}\cdot\vec{z}) \\
& = 
(\vec{v}\cdot\vec{u})\vec{w}\wedge\vec{z} - 
\vec{u}\wedge(\vec{v}\cdot(\vec{w}\wedge
\vec{z})) . 
\end{split}
\end{equation}
So, we have defined a product $\vec{v}\cdot A_k \in \bigwedge^{k-1}(\mathbb{R}^3)$ for $A_k \in \bigwedge^k(\mathbb{R}^3)$, which we call the {\em interior product}. This product is extended to all $\bigwedge(\mathbb{R}^3)$ by linearity. We can also generalize this product by
\begin{equation}
\label{def.int.biv}
(\vec{v}\wedge\vec{u})\cdot A_k = 
\vec{v}\cdot(\vec{u}\cdot A_k) , 
\end{equation}
when $k \geq 2$, and 
\begin{equation}
(\vec{v}\wedge\vec{u}\wedge\vec{w})\cdot A_k = 
\vec{v}\cdot(\vec{u}\cdot(\vec{w}\cdot A_k)) , 
\end{equation}
for $k = 3$. These expressions define the interior product 
$A_k \cdot B_j$ for $k \leq j$. For $k > j$, we define 
\begin{equation}
A_k \cdot B_j = (-1)^{j(k-1)} B_j \cdot A_k . 
\end{equation}
For more details, see~\cite{VazRocha,DFM,lounesto}. 

We define the scalar product $(A_k|B_j)$ as 
\begin{equation}
(A_k | B_j) = \begin{cases} 0 , & \quad k \neq j , \\
\tilde{A}_k \cdot B_j , & \quad k = j =1,2,3,\\
A_k  B_j , & \quad k = j =0.
\end{cases}
\end{equation}

\subsection{The Hodge Star Operator}

The dimension of $\bigwedge^k(\mathbb{R}^3)$ is $\binom{3}{k}$. Since
$\binom{n}{k} = \binom{n}{n-k}$, then 
$\bigwedge^k(\mathbb{R}^3)$ and $\bigwedge^{3-k}(\mathbb{R}^3)$ have
the same dimension, that is, they are isomorphic vector spaces. 
There is, however, no canonical isomorphism between these spaces, 
which means that one such isomorphism has to be defined. An
important one is the {\em Hodge star} isomorphism $\star: 
\bigwedge^k(\mathbb{R}^3) \rightarrow \bigwedge^{n-k}(\mathbb{R}^3)$ 
defined as 
\begin{equation}
\label{hodge.def}
A \wedge \star B = (A|B) \Omega, \quad \forall A \in \bigwedge(\mathbb{R}^3) , 
\end{equation}
where 
\begin{equation}
\Omega = \vec{e}_1\wedge \vec{e}_2 \wedge \vec{e}_3. 
\end{equation}
It can be proven~\cite{VazRocha} that this definition is equivalent to 
\begin{equation}
\label{hodge.equiv}
\star A_k = \tilde{A}_k \cdot \Omega , \qquad \star 1 = \Omega, 
\end{equation}
for $A_k \in \bigwedge^k(\mathbb{R}^3)$. 
It follows then that 
\begin{xalignat}{2}
\label{hodge.1}
& \star 1 = \vec{e}_1 \wedge \vec{e}_2\wedge \vec{e}_3 , 
& \quad & \star \vec{e}_1 = \vec{e}_2 \wedge\vec{e}_3 , \\
\label{hodge.2}
& \star \vec{e}_2 = \vec{e}_3\wedge\vec{e}_1, & \quad & 
\star \vec{e}_3 = \vec{e}_1\wedge\vec{e}_2 , \\
\label{hodge.3}
& \star \vec{e}_1\wedge\vec{e}_2 = \vec{e}_3 , & \quad 
& \star \vec{e}_3\wedge\vec{e}_1 = \vec{e}_2 , \\
\label{hodge.4}
& \star \vec{e}_2\wedge\vec{e}_3 = \vec{e}_1 , & \quad 
& \star \vec{e}_1\wedge\vec{e}_2 \wedge\vec{e}_3 = 1 . 
\end{xalignat}

\section{Paravectors}
\label{sec.3}

Let us fix the notation $\bigwedge^{-1}(\mathbb{R}^3) = \emptyset$ and $\bigwedge^{4}(\mathbb{R}^3) = \emptyset$. 
We define $\prod^k(V)$ as 
\begin{equation}
\textstyle{\prod^k}(\mathbb{R}^3) = \textstyle{\bigwedge^{k-1}}(\mathbb{R}^3) \oplus 
\textstyle{\bigwedge^k}(\mathbb{R}^3) , \qquad 
k = 0, 1, 2, 3, 4. 
\end{equation}  
We call the elements of $\prod^k(V)$ \mbox{$k$-paravectors} for $k = 1, 2, 3, 4$, where for $k=0$ paravectors are scalars. A $1$-paravector is called 
simply a paravector, and therefore a paravector is the sum of a $0$-vector (scalar) 
and a vector; a  biparavector is a sum of a vector and a bivector; 
a triparavector is a sum of a bivector and a trivector, and a 
quadriparavector is just a trivector. 
To establish a notation that clearly distinguishes $k$-paravectors and $k$-vectors, we will denote elements of $\bigwedge^k(\mathbb{R}^3)$ by 
capital Roman letters with sub-index $k$, like $A_k$, and we will denote elements of 
$\prod^k(V)$ by capital Roman letters in sans serif fonts with sub-index 
between curly brackets, like 
$\mathsf{A}_{\{k\}}$. Using this notation, an arbitrary $k$-paravector is an 
element of the form 
\begin{equation}
\mathsf{A}_{\{k\}} = A_{k-1} + A_k , \qquad k = 0, 1, 2, 3, 4, 
\end{equation}
where $A_{-1} = 0$ and $A_{n+1} = 0$. 

Given a $k$-paravector, we can extract its $(k-1)$-vector and
$k$-vector parts using grade involution. In fact, we have
\begin{gather}
 \langle \mathsf{A}_{\{k\}}\rangle_{k-1} = \frac{1}{2}\left[\mathsf{A}_{\{k\}} - 
(-1)^{k}  \widehat{\mathsf{A}}_{\{k\}}\right] , \\
 \langle \mathsf{A}_{\{k\}}\rangle_{k} = \frac{1}{2}\left[\mathsf{A}_{\{k\}} + 
(-1)^{k} \widehat{\mathsf{A}}_{\{k\}}\right] .
\end{gather}

In three dimensions, $k$-paravectors ($k=0,1,2,3,4$) can also be  defined using 
the reversion and conjugation operations. Indeed, given an arbitrary 
$\phi \in \bigwedge(\mathbb{R}^3)$, we have 
\begin{equation}
\begin{cases}
 \tilde{\phi} = \phi & \Rightarrow \phi \text{ is a paravector} , \\
 \bar{\phi} = -\phi & \Rightarrow \phi \text{ is a biparavector} , \\
 \tilde{\phi} = -\phi & \Rightarrow \phi \text{ is a triparavector} , \\
 \bar{\phi} = \phi & \Rightarrow 
\phi \text{ is a sum of 0-paravector and 4-paravector} .  
\end{cases}
\end{equation}
In the last case, when $\overline{\phi} = \phi$, 
 if also $\hat{\phi} = \phi$, then $\phi$ is a scalar, 
 while if $\hat{\phi} = - \phi$, then $\phi$ is a 
 quadriparavector.

Just like we have an 
exterior product $\wedge$ that gives a $(k+l)$-vector from 
the product of a $k$-vector and a $l$-vector, we would like 
to have a product that gives a $(k+l)$-paravector from the
product of a $k$-paravector and a $l$-paravector. We denote this product
by $\curlywedge$. 
It is natural to suppose that $\curlywedge$ could be written in terms 
of the exterior product $\wedge$. In fact, note that 
\begin{equation}
\begin{split}
\mathsf{A}_{\{k\}}\wedge \mathsf{B}_{\{l\}} & = 
(A_{k-1}+A_k)\wedge(B_{l-1}+B_l) \\
& = 
A_{k-1}\wedge B_{l-1} + A_{k-1}\wedge B_l + A_k \wedge B_{l-1} + 
A_k \wedge B_l ,
\end{split}
\end{equation} 
where 
\begin{alignat}{2}
& A_{k-1}\wedge B_{l-1} \in \textstyle{\bigwedge^{k+l-2}}(\mathbb{R}^3) , & \quad 
& A_{k-1}\wedge B_l \in \textstyle{\bigwedge^{k+l-1}}(\mathbb{R}^3) , \\ 
& A_k \wedge B_{l-1} \in \textstyle{\bigwedge^{k+l-1}}(\mathbb{R}^3) , & \quad 
& A_k \wedge B_l \in \textstyle{\bigwedge^{k+l}}(\mathbb{R}^3) . 
\end{alignat}
We can see that there is a term, namely $A_{k-1}\wedge B_{l-1}$, that
does not belong to $\prod^{k+l}(\mathbb{R}^3)$. One way to fix this problem
is to use projectors.  Denote 
\begin{equation}
\langle \quad \rangle_{\{k\}} = \langle \quad \rangle_{k-1} + 
\langle \quad \rangle_k . 
\end{equation} 
Then 
\begin{equation}
\langle \mathsf{A}_{\{k\}}\wedge \mathsf{B}_{\{l\}} \rangle_{\{k+l\}} = 
A_{k-1}\wedge B_l + A_k \wedge B_{l-1} + 
A_k \wedge B_l ,
\end{equation}
which is a $(k+l)$-paravector. Note also that, because of the associativity 
of the exterior product, we have  
\begin{equation}
\begin{split}
\left\langle \langle \mathsf{A}_{\{k\}}\wedge \mathsf{B}_{\{l\}}\rangle_{\{k+l\}} 
\wedge \mathsf{C}_{\{m\}} \right\rangle_{\{k+l+m\}} & = 
\left\langle \mathsf{A}_{\{k\}}\wedge \langle \mathsf{B}_{\{l\}} 
\wedge \mathsf{C}_{\{m\}}\rangle_{\{l+m\}} \right\rangle_{\{k+l+m\}} \\
& = A_{k-1}\wedge B_l \wedge C_m + A_k \wedge B_{l-1} \wedge C_m \\
& \quad + A_k \wedge B_l \wedge C_{m-1} + A_k \wedge B_l \wedge C_m .
\end{split}
\end{equation}
These results suggest defining the (associative) exterior product 
of paravectors $\curlywedge$ as
\begin{equation}
\mathsf{A}_{\{k\}}\curlywedge \mathsf{B}_{\{l\}} = 
\langle \mathsf{A}_{\{k\}}\wedge \mathsf{B}_{\{l\}}\rangle_{\{k+l\}} . 
\end{equation}

To interpret the result of the product of paravectors, we 
need first an interpretation for a paravector. 
Consider a paravector 
$\mathsf{P}$ of the form 
\begin{equation}
\mathsf{P} = 1 + \vec{p} .
\end{equation}
We interpret this paravector as describing a point $P$ in an affine 
space with coordinates in relation to the origin (described by $\mathsf{E} = 1$)  
given by the coordinates of the vector $\vec{p}$. 
An arbitrary paravector of the form 
\begin{equation}
\mathsf{X} = x_0 + \vec{x} 
\end{equation}
is interpreted as a weighted point, with weight $x_0$ and
located at $\vec{x}/|x_0|$ (the reason for the
use of the absolute value will be clear below). Then the points 
$\mathsf{P}$ and $m\mathsf{P}$
have the same location, but different weights $1$ and $m$, 
respectively. Some authors call these object points with mass
or massive points, but since $m$ can be negative, we prefer 
to use the terminology {\em weight}.  We will say that 
points with $x_0 > 0$ have positive orientation, 
while points with $x_0 < 0$ have negative orientation. 

The sum of a point $\mathsf{P}$ and a vector $\vec{v}$ 
gives another point $\mathsf{Q}$ located at $\vec{q} = 
\vec{p} + \vec{v}$, 
\begin{equation}
\mathsf{Q} = \mathsf{P} + \vec{v} = 
1 + \vec{p} + \vec{v} = 1 + \vec{q} . 
\end{equation}
The sum of a point $m_1\mathsf{P}_1$ with weight $m_1$ and 
a point $m_2\mathsf{P}_2$ with weight $m_2$ is 
\begin{equation}
m_1\mathsf{P}_1 + m_2\mathsf{P}_2 = (m_1+m_2) + (m_1 \vec{p}_1 + 
m_2 \vec{p}_2)  ,
\end{equation}
which is interpreted as a 
 point with weight $m_1 + m_2$ located at the center of 
mass $\vec{p}_{\scriptstyle \text{CM}} = 
(m_1 \vec{p}_1 + 
m_2 \vec{p}_2)  /(m_1+ m_2)$. When the weights have 
opposite signs, that is, $m_1 = - m_2$, the result is a vector 
$m_1(\vec{p}_1 - \vec{p}_2) = 
m_2(\vec{p}_2 - \vec{p}_1)$.

Now let us consider the product of paravectors. 
We already know that $\mathsf{P}\curlywedge\mathsf{Q}$ is 
a biparavector, but we would like to have an interpretation 
of a product of paravectors as a kind of product of points. 
We would like to interpret 
the result of this product as representing a line segment with 
orientation. Given points $P$ and $Q$, we would like
to distinguish the line segment leaving the point $P$ and
reaching the point $Q$ from the line segment leaving the
point $Q$ and reaching the point $P$. Moreover, in the 
limit where $Q$ approaches $P$, this product must 
be null. If we write $\mathsf{P} = p_0 + \vec{p}$ and
$\mathsf{Q} = q_0 + \vec{q}$, the condition 
$\mathsf{P}\curlywedge\mathsf{Q} = 0$ implies that 
$\vec{q} = \alpha \vec{p}$ and $q_0 = -\alpha p_0$, 
for $\alpha \in \mathbb{R}$, that is, 
$\mathsf{Q}$ must be of the form $\alpha(-p_0 + \vec{p})$. 
According to our interpretation, 
the paravectors $p_0 + \vec{p}$ and $-p_0 + \vec{p}$ represent 
points with opposite orientations. Therefore, to 
have a product of paravectors that vanishes 
when the paravectors are the same, the change
of orientation has to be taken into account in our product, 
and a useful way of doing this is by using an 
algebraic operation. 

Consider the points represented by the paravectors 
$\mathsf{P} = 1 + \vec{p}$ and $\mathsf{Q} = 1 + 
\vec{q}$. The points with opposite weights are
given by $\mathsf{P}^\dagger = -1 + \vec{p}$ and 
$\mathsf{Q}^\dagger = -1 + \vec{q}$. It would be 
useful if $\mathsf{P}^\dagger$  and $\mathsf{Q}^\dagger$ could 
be obtained from $\mathsf{P}$ and $\mathsf{Q}$ by means of
some of the algebraic involutions already discussed, since 
the change of orientation is obviously an operation whose
composition with itself is the identity operation. One such 
operation is related to conjugation, that is, $\bar{\mathsf{P}} = 
1 - \vec{p}$ and $\bar{\mathsf{Q}} = 1 - \vec{q}$. 
Then $\mathsf{P}^\dagger = - \bar{\mathsf{P}}$ and 
$\mathsf{Q}^\dagger = -\bar{\mathsf{Q}}$. Note that
the operation $\mathsf{P} \mapsto -\mathsf{P}$ changes
not only the orientation of the point but also its location, from
$\vec{p}$ to $-\vec{p}$. The use of 
conjugation restores the point to its original location, 
changing only its orientation. For the sake of
convenience, we define 
\begin{equation}
(\mathsf{P}_1\curlywedge\cdots\curlywedge\mathsf{P}_n)^\dagger = 
\mathsf{P}_1^\dagger \curlywedge\cdots \curlywedge \mathsf{P}_1^\dagger, 
\qquad \mathsf{P}^\dagger = -\bar{\mathsf{P}} . 
\end{equation}

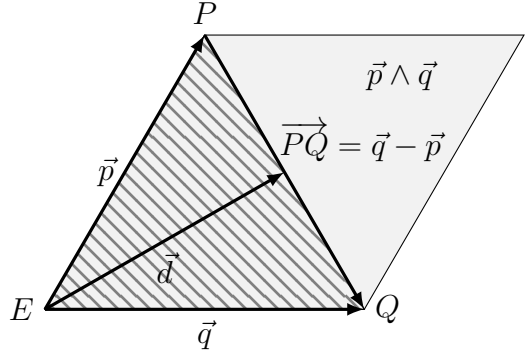
\begin{figure}
\begin{center}
\begin{tikzpicture}[>=latex,xscale = 2.1,yscale=2.1]
\fill [gray!10]
      (0, 0) -- (2,0) -- (3,1.732) -- (1,1.732) -- cycle;
\draw (0, 0) -- (2,0) -- (3,1.732) -- (1,1.732) -- cycle;
\draw [pattern= custom north west lines,hatchspread=6pt,hatchthickness=1pt,hatchcolor= gray] (0, 0) -- (2,0) -- (1,1.732) -- cycle;
\draw [very thick,->] (0,0) -- (2,0);
\draw [very thick,->] (0,0) -- (1,1.732);
\draw [very thick,->] (1,1.732) -- (2,0);
\draw [very thick,->] (0,0) -- (1.5,0.866);
\node [below] at (1,0) {$\vec{q}$};
\node [left] at (0.5,0.866) {$\vec{p}$};
\node [below] at (0.75,0.43) {$\vec{d}$};
\node [above right] at (1.4,0.866) {$\overrightarrow{PQ} = 
\vec{q}-\vec{p}$};
\node [right] at (2,0) {$Q$};
\node [above] at (1,1.732) {$P$};
\node [left] at (0,0) {$E$};
\node [below left] at (2.5,1.632) {$\vec{p}\wedge\vec{q}$};
\end{tikzpicture}
\end{center}
\caption{The information contained in the 
biparavector representation $\mathcal{L}$ 
of the line segment 
$PQ$: the parallelogram described by the bivector $\vec{p}\wedge
\vec{q}$, the oriented line segment $\vec{q}-\vec{p}$ 
and the support vector $\vec{d}$. } 
\label{fig:biparavectorline}
\end{figure}

Let us consider therefore the product $\mathsf{P}\curlywedge\mathsf{Q}^\dagger$, 
which gives 
\begin{equation}
\label{eq.line}
\mathcal{L} = \mathsf{P}\curlywedge\mathsf{Q}^\dagger = \vec{q} - 
\vec{p} + \vec{p}\wedge\vec{q} .
\end{equation}
The vector part of this biparavector is $\langle \mathcal{L}\rangle_1 = \overrightarrow{PQ} = 
\vec{q} - \vec{p}$, 
and its bivector part is $\langle \mathcal{L}\rangle_2 = 
M = \vec{p}\wedge\vec{q}$, 
which is interpreted as the moment of the line about the origin. 
This is the biparavector representation of the Pl\"ucker 
coordinates of a line~\cite{DFM}. In fact, the coordinates
of $M = \vec{p}\wedge\vec{q}$ are the same as the 
coordinates of $\vec{p}\times\vec{q}$. Moreover, 
$|M| = 2 A$, where $A$ is the area of the
triangle $EPQ$, where $E$ is the origin. So, $|M| = 2 \frac{1}{2} 
|\vec{d}| |\overrightarrow{PQ}|$, 
where $|\vec{d}|$ is the distance from the line segment 
$PQ$ to the origin when 
$\vec{d}\cdot \overrightarrow{PQ} = 0$. 
This vector (the support vector) is $\vec{d} = M \cdot \overrightarrow{PQ} / 
|\overrightarrow{PQ}|^2$ (the fact that 
in this case $\vec{d}\cdot \overrightarrow{PQ} = 0$
can be seen from eq.\eqref{def.int.biv}, for example), and 
$M = \overrightarrow{PQ} \wedge \vec{d}$ 
is the counterpart of the definition of the moment
of the line as $\overrightarrow{PQ}\times \vec{d}$. 
It is also convenient to write the support vector $\vec{d}$ as  
\begin{equation}
\vec{d} = 
\frac{\langle \mathcal{L}\rangle_2 \cdot \langle \mathcal{L}\rangle_1}
{|\langle \mathcal{L}\rangle_1|^2} . 
\end{equation}
See Figure~\ref{fig:biparavectorline}.

One can think of defining a biparavector $\mathsf{P}^\dagger \curlywedge\mathsf{Q}$ 
and interpret this paravector as describing a line segment ending in $P$ and 
starting at $Q$. The product $\mathsf{P}^\dagger \curlywedge \mathsf{Q}$ 
results in $\overrightarrow{QP} + M$, where $\overrightarrow{QP} = 
\vec{p} - \vec{q}$ and $M = \vec{p}\wedge\vec{q}$. 
Although the presence of the vector $\overrightarrow{QP}$ is in 
accordance with our tentative interpretation, there is a problem 
with the interpretation of the moment of the line. If we define, 
as discussed above, the vector $\vec{d}_1 = 
M \cdot \overrightarrow{QP} /|\overrightarrow{QP}|^2$, 
we obtain $\vec{d}_1 = -\vec{d}$, which 
is not what we expect, which would be $\vec{d}_1 = \vec{d}$. 
One way of fixing this problem is to define, in this case, $\vec{d}^\prime = 
 \overrightarrow{QP} \cdot M /|\overrightarrow{QP}|^2$, and 
then $\vec{d}^\prime = \vec{d}$. 
The question here is one of consistency: we can choose to represent
a line segment by biparavectors of the form $\mathsf{P} \curlywedge \mathsf{Q}^\dagger$ 
or $\mathsf{P}^\dagger \curlywedge \mathsf{Q}$, but if we want to 
use both at the same time, we need to take some care. 
If we choose to represent an oriented line segment by a product 
of points represented by paravectors as in $\mathsf{P} \curlywedge \mathsf{Q}^\dagger$, 
the best thing, in our opinion, is to continue with this interpretation and
not mix the different choices. In this case, 
the biparavector represented starting in $Q$ and ending in 
$P$ is $\mathsf{Q} \curlywedge \mathsf{P}^\dagger = 
\overrightarrow{QP} - M$, and it is such 
that 
\begin{equation}
\mathsf{Q} \curlywedge \mathsf{P}^\dagger = 
(\mathsf{P} \curlywedge \mathsf{Q}^\dagger)^\dagger ,
\end{equation}
which is consistent with our interpretation of $\dagger$ as related
to orientation. Note that 
\begin{equation}
\label{support.1}
\mathsf{P}^\dagger\curlywedge\mathsf{Q} = 
\widetilde{(\mathsf{Q} \curlywedge \mathsf{P}^\dagger)} .
\end{equation}

Now consider the product of three paravectors, 
$\mathsf{P}\curlywedge\mathsf{Q}^\dagger
\curlywedge\mathsf{R}$, which results in 
\begin{equation}
\mathsf{P}\curlywedge\mathsf{Q}^\dagger
\curlywedge\mathsf{R} = \vec{p}\wedge\vec{q} 
- \vec{p}\wedge\vec{r} + \vec{q}\wedge
\vec{r} + \vec{p}\wedge \vec{q}\wedge
\vec{r} . 
\end{equation} 
In contrast to the case $\mathsf{P}\curlywedge\mathsf{Q}^\dagger$, which
is always non-null if the points are different, we can have 
the situation where $\mathsf{P}\curlywedge\mathsf{Q}^\dagger
\curlywedge\mathsf{X} = 0$. For this to happen, we must have 
\begin{gather}
\vec{p}\wedge \vec{q}\wedge
\vec{x} = 0 , \\
\vec{p}\wedge\vec{q} 
- \vec{p}\wedge\vec{x} + \vec{q}\wedge
\vec{x} = 0 . 
\end{gather}
From the first equation, we conclude that
\begin{equation}
\vec{x} = s \vec{p} + t \vec{q} , 
\end{equation}
where $t$ and $s$ are scalars, and using this in the 
second equation, we obtain
\begin{equation}
(1 - s - t) \vec{p} \wedge\vec{q} = 0 , 
\end{equation}
from which we conclude that 
\begin{equation}
s + t = 1 , 
\end{equation}
and then 
\begin{equation}
\vec{x} = \vec{p} + t(\vec{q}-\vec{p}) , 
\end{equation}
where $t$ is a scalar. Then $\mathsf{X} = 1 + \vec{x}$ is a 
point along the line passing through the points $P$ and $Q$. 
In summary: the biparavector $\mathcal{L} = \mathsf{P}\curlywedge
\mathsf{Q}^\dagger$ describes the line segment from $P$ to $Q$ and
the equation of the line that passes through these points is 
\begin{equation}
\mathcal{L}\curlywedge \mathsf{X} = 0 . 
\end{equation}

\begin{figure}
\begin{center}
\begin{tikzpicture}[>=latex,xscale = 1.5,yscale=1.5]
\fill [gray!10]     (0, 0) -- (2,0) -- (2.5,0.866) -- (2.5,3.866) --
(0.5,3.866) -- (0,3) -- cycle;
\draw (0, 0) -- (2,0) -- (2,3) -- (0,3) -- cycle;
\draw (2+0.5,0+0.866) -- (2+0.5,3+0.866) -- (0+0.5,3+0.866);
\draw [dashed] (0.5,3.866) -- (0+0.5, 0+0.866) -- (2.5,0.866);
\draw [pattern= custom north west lines,hatchspread=6pt,hatchthickness=1pt,hatchcolor= gray] (0.5,0.866) -- (0,3) -- (2,0) -- cycle;
\draw [very thick,->] (0,0) -- (2,0);
\draw [very thick,->] (0,0) -- (0,3);
\draw [very thick,->] (0,0) -- (0.5,0.866);
\draw (2,0) -- (2.5,0.866);
\draw (2,3) -- (2.5,3.866);
\draw (0,3) -- (0.5,3.866);
\draw [very thick, ->] (0,0) -- (1,0.866);
\draw [very thick] (0.5,0.866) -- (0,3) -- (2,0) -- cycle;
\node [below] at (1,0) {$\vec{r}$};
\node [left] at (0,1.5) {$\vec{p}$};
\node [above] at (0.15,0.38) {$\vec{q}$};
\node [below] at (0.5,0.43) {$\vec{d}$};
\node [right] at (2,0) {$R$};
\node [left] at (0,3) {$P$};
\node [left] at (0.5,0.866) {$Q$};
\node [left] at (0,0) {$E$};
\node [above] at (1.2,2) {$\langle \mathcal{P}\rangle_2$};
\draw [dashed, ->] (1.2,2) to [bend left] (0.8,1.3);
\node [right] at (0.6,3.4) {$\vec{p}\wedge\vec{q}\wedge\vec{r}$};
\end{tikzpicture}
\end{center}
\caption{The information contained in the 
triparavector representation $\mathcal{P}$ 
of the plane fragment 
$PQR$: the parallelepiped described by the trivector $\vec{p}\wedge
\vec{q}\wedge\vec{r}$, the oriented plane 
fragment described by the bivector $\langle \mathcal{P}\rangle_2 = 
\vec{p}\wedge\vec{q} 
- \vec{p}\wedge\vec{r} + \vec{q}\wedge
\vec{r} $ 
and the support vector $\vec{d}$. } 
\label{fig:triparavectorplane}
\end{figure}
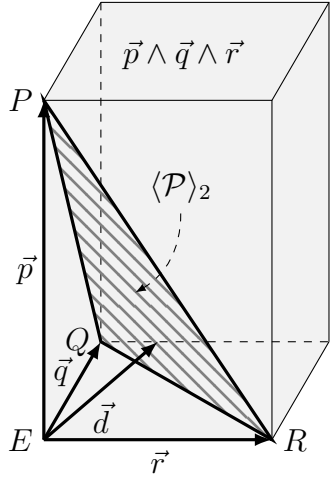

When $\mathcal{L}\curlywedge \mathsf{X}$ is non-null, the product 
\begin{equation}
\label{plane.1}
\mathcal{P} = \mathsf{P}\curlywedge\mathsf{Q}^\dagger
\curlywedge\mathsf{R} 
\end{equation}
with 
\begin{gather}
\label{plane.2}
\langle \mathcal{P}\rangle_2 =  \vec{p}\wedge\vec{q} 
- \vec{p}\wedge\vec{r} + \vec{q}\wedge
\vec{r}, \\
\label{plane.3}
\langle \mathcal{P}\rangle_3 =  \vec{p}\wedge \vec{q}\wedge
\vec{r},  
\end{gather}
is the triparavector representation of the plane fragment defined by 
the points $P$, $Q$ and $R$, just like $\mathcal{L}$ is the biparavector 
representation of the line segment defined by $P$ and $Q$. The bivector  
$\langle \mathcal{P}\rangle_2$ describes the direction and the
orientation of the plane, and the trivector $\langle \mathcal{P}\rangle_3$ 
is a kind of moment of the plane about the origin. Its absolute value
$|\langle \mathcal{P}\rangle_3|$ is a measure of the distance from 
the plane to the origin. In fact, $|\langle \mathcal{P}\rangle_3|$ is 
the volume of the parallelepiped defined by the vectors 
$\vec{p}$, $\vec{q}$ and $\vec{r}$, which 
is $6$ times the volume $V$ of the tetrahedron defined by the points
$EPQR$, where $E$ is the origin. However, $V = A d/3$, where $A$ is the
area of the base (the triangle $PQR$) and $d$ is the height of the
tetrahedron (the distance from the plane fragment to the origin). 
Since $A = |\langle \mathcal{P}\rangle_2|/2$, we have 
$|\langle \mathcal{P}\rangle_3| = |\langle \mathcal{P}\rangle_2| d$. 
The support vector $\vec{d}$ is such that $|\vec{d}\,| = d$ 
and is orthogonal to the plane fragment, so we have 
\begin{equation}
\vec{d} = \frac{\langle \mathcal{P}\rangle_3 \cdot 
\langle \mathcal{P}_2\rangle}{|\langle \mathcal{P}\rangle_2|^2} ,
\end{equation}
which is to be compared with eq.\eqref{support.1}; see Figure~\ref{fig:triparavectorplane}. 

Just like $\mathcal{L}\curlywedge \mathsf{X} = 0$ is the
equation of the line that passes through $P$ and $Q$, we 
expect that $\mathcal{P}\curlywedge\mathsf{X}^\dagger = 0$ 
to be the equation of the plane that passes through the non-collinear 
points $P$, $Q$ and $R$. In fact, in this case we have 
\begin{gather}
\vec{q}\wedge\vec{r}\wedge\vec{x} - 
\vec{p}\wedge\vec{r}\wedge\vec{x} + 
\vec{p}\wedge\vec{q}\wedge\vec{x} - 
\vec{p}\wedge\vec{q}\wedge\vec{r} = 0 , \\
\vec{p}\wedge\vec{q}\wedge\vec{r}\wedge
\vec{x} = 0 . 
\end{gather} 
The last condition is trivial in a three dimensional space. 
Let us suppose that $\vec{p}$, $\vec{q}$ and $\vec{r}$ 
are linearly independent. Then 
\begin{equation}
\vec{x} = s\vec{p} + t\vec{q} + u \vec{r} , 
\end{equation}
and using this on the second last condition, we have 
\begin{equation}
(s + t + u - 1) \vec{p}\wedge\vec{q}\wedge\vec{r} = 0 , 
\end{equation}
that is, 
\begin{equation}
\vec{x} = \vec{p} + t(\vec{q}-\vec{p}) + 
u (\vec{r} - \vec{p}) , 
\end{equation}
which is the vector form of the equation of a plane. 
Then $\mathsf{X} = 1 + \vec{x}$ is a 
point in the plane passing through the points $P$, $Q$ and $R$.

Due to the associativity of the product $\curlywedge$ of paravectors, 
we can consider the product $\mathsf{P}\curlywedge 
\mathsf{Q}^\dagger\curlywedge\mathsf{R}\curlywedge\mathsf{S}^\dagger$ as 
the product of two biparavectors $\mathcal{L}$ and $\mathcal{M}$ 
representing two  lines, 
given by 
\begin{gather} 
\mathcal{L} = \mathsf{P}\curlywedge\mathsf{Q}^\dagger = 
\overrightarrow{PQ} + \vec{p}\wedge\vec{q} , \\
\mathcal{M} = \mathsf{R}\curlywedge\mathsf{S}^\dagger = 
\overrightarrow{RS} + \vec{r}\wedge\vec{s}. 
\end{gather}
Then $\mathcal{L}\curlywedge\mathcal{M} = 0$ means that 
the two lines lie in the same plane. Since 
\begin{equation}
\vec{q}\wedge\vec{r}\wedge\vec{s} - 
\vec{p}\wedge\vec{r}\wedge\vec{s} + 
\vec{p}\wedge\vec{q}\wedge\vec{s} - 
\vec{p}\wedge\vec{q}\wedge\vec{r} = 
-\overrightarrow{PR} \wedge \overrightarrow{PQ}\wedge
\overrightarrow{RS} , 
\end{equation}
the condition $\mathcal{L}\curlywedge\mathcal{M} = 0$ 
implies that 
\begin{equation}
\overrightarrow{PR} \wedge \overrightarrow{PQ}\wedge
\overrightarrow{RS} = 0 . 
\end{equation}
This condition is satisfied if any of the products 
$\overrightarrow{PQ}\wedge \overrightarrow{RS}$ , 
$\overrightarrow{PR}\wedge\overrightarrow{PQ}$ or
$\overrightarrow{PR}\wedge\overrightarrow{RS}$ vanishes. 
Since $\overrightarrow{PQ}$ and $\overrightarrow{RS}$ are
the vectors that define the directions of the lines 
$\mathcal{L}$ and $\mathcal{M}$, respectively, when 
these lines lie in the same plane we have 
following possibilities:
\begin{enumerate}
  \item[(i)] if 
$\overrightarrow{PQ}\wedge \overrightarrow{RS} = 0$ and 
$\overrightarrow{PR}\wedge\overrightarrow{PQ} \neq 0$ and
$\overrightarrow{PR}\wedge\overrightarrow{RS} \neq 0$, then $\mathcal{L}$ and $\mathcal{M}$
    are parallel lines;
    \item [(ii)] $\overrightarrow{PQ}\wedge \overrightarrow{RS} = 0$ and 
$\overrightarrow{PR}\wedge\overrightarrow{PQ} = 0$ and
$\overrightarrow{PR}\wedge\overrightarrow{RS} = 0$, then $\mathcal{L}$ and $\mathcal{M}$  
      are coincident lines;
      \item [(iii)] if $\overrightarrow{PQ}\wedge \overrightarrow{RS} \neq 0$, then 
$\mathcal{L}$ and $\mathcal{M}$ are intersecting lines---and if $\overrightarrow{PQ}\cdot \overrightarrow{RS} = 0$ then
$\mathcal{L}$ and $\mathcal{M}$ are perpendicular lines.
\end{enumerate}

The quadriparavector $\mathcal{V} = \mathsf{P}\curlywedge 
\mathsf{Q}^\dagger\curlywedge\mathsf{R}\curlywedge\mathsf{S}^\dagger \neq 0$ 
is 
\begin{equation}
\begin{split}
\mathcal{V} & = \mathsf{P}\curlywedge 
\mathsf{Q}^\dagger\curlywedge\mathsf{R}\curlywedge\mathsf{S}^\dagger = 
\vec{q}\wedge\vec{r}\wedge\vec{s} - 
\vec{p}\wedge\vec{r}\wedge\vec{s} + 
\vec{p}\wedge\vec{q}\wedge\vec{s} - 
\vec{p}\wedge\vec{q}\wedge\vec{r} \\
& = (\vec{q}-\vec{p})\wedge
(\vec{r}-\vec{q})\wedge
(\vec{s}-\vec{r}) = 
\begin{vmatrix} p^1 & p^2 & p^3 & 1 \\ q^1 & q^2 & q^3 & 1 \\
r^1 & r^2 & r^3 & 1 \\ s^1 & s^2 & s^3 & 1 \end{vmatrix}
\vec{e}_1\wedge\vec{e}_2\wedge\vec{e}_3 . 
\end{split}
\end{equation}
This determinant is the volume of the parallelepiped 
defined, for example,  by the vectors $\overrightarrow{PQ} = \vec{q}-\vec{p}$, 
$\overrightarrow{QR} = \vec{r}-\vec{q}$ and 
$\overrightarrow{RS} = \vec{s}-\vec{r}$, or 
$6$ times the volume of the tetrahedron defined by the points 
$P$, $Q$, $R$ and $S$. There is no $k$-paravector with 
$k > 4$ in a three dimensional vector space.

Finally, we observe that the Hodge dual operator can have
an interesting role in the representation of geometric 
objects using paravectors. Let us consider a line 
represented by $\mathcal{L} = \vec{l} + M$---see eq.\eqref{eq.line}. 
The dual of the bivector $M$ is a vector $\vec{m}$. Then 
we can write
\begin{equation}
\mathcal{L} = \vec{l} + \star\vec{m} , 
\end{equation}
which is the usual Pl\"ucker representation of a line, 
with homogeneous coordinates $[\vec{l},\vec{m}]$, although
with a different geometric interpretation. In the same 
manner, let us consider a plane represented by 
$\mathcal{P}$. Using the Hodge dual operator, we can
write the bivector part $\langle \mathcal{P}\rangle_2$ 
as $\star \vec{n}$ and the trivector part as $c \star 1$, 
that is, 
\begin{equation}
\label{plane.dual}
\mathcal{P} = \star(\vec{n} + c) . 
\end{equation}
The vector $n$ is interpreted as a normal vector to 
the plane and $|c|$ is interpreted as the volume 
of a parallelepiped defined by the vectors constructed 
from three points on the plane
and the origin. The equation of the plane $\mathcal{P}\curlywedge
\mathsf{X}^\dagger = 0$ can be written therefore as 
\begin{equation}
-(\star\vec{n})\wedge \vec{x} + c\star 1 = 0 , 
\end{equation}
and from eq.\eqref{hodge.def}, 
\begin{equation}
\label{eq.plane.usual}
-(\vec{n}\cdot\vec{x})\star 1 + c \star 1 = 0 , 
\end{equation}
which gives the equation of the plane in the usual 
form $\vec{n}\cdot\vec{x} = c$.

We summarize these ideas in Table~\ref{tab:summary}

\begin{table}[h]
\begin{center}
\begin{tabular}{| >{\centering\arraybackslash} m{4cm}  | >{\centering\arraybackslash} m{4cm}  | >{\centering\arraybackslash} m{5cm}  |}
\hline \hline  
Geometric object & Algebraic object & Multivector expression \\ 
\hline \hline 
Point & Paravector 
$\mathsf{P}$ & $\mathsf{P} = 1 + \vec{p}$ \\ 
\hline 
 Line segment & Biparavector $\mathcal{L} = 
\mathsf{P}\curlywedge \mathsf{Q}^\dagger$ & 
$ \mathcal{L} \!\begin{aligned}[t]  & = \vec{q}-\vec{p} + \vec{p}\wedge \vec{q} \\
&  = \vec{l} + \star \vec{m} \end{aligned} $ \\
 \hline 
 Plane fragment & Triparavector $\mathcal{P} = 
\mathsf{P}\curlywedge \mathsf{Q}^\dagger \curlywedge 
\mathsf{R}$ & $\mathcal{P} \!\begin{aligned}[t]
    &= \vec{p}\wedge\vec{q} - \vec{p}\wedge\vec{r} + \vec{q}\wedge\vec{r}\\
    & \qquad + \vec{p}\wedge\vec{q}\wedge\vec{r}  \\
    & = \star \vec{n} + c \star 1 
    \end{aligned}$ \\
\hline \hline 
\end{tabular}
\end{center}
\caption{Summary of the relation between geometric objects and $k$-paravectors.}
\label{tab:summary}
\end{table}

\vfill 

\section{The Algebra of Transformations}
\label{sec.4}

Our objective now is to study the transformations of points, lines and
planes in terms of their representation using paravectors, 
biparavectors and triparavectors. We can do this through 
combinations of the exterior and interior products, but 
there is a problem: the exterior algebra, which is at 
the base of our description, is not a matrix algebra. 
This matter is not a problem in itself, but from a practical 
point of view, it would be better to work with a 
matrix algebra. There is also the annoying 
non-associativity of exterior and interior products; 
for example, suppose we want to perform an interior  
product from the left by $\vec{e}_2$ and 
an exterior product from the right by $\vec{e}_1\wedge\vec{e}_2$. 
It is easy to see that
$(\vec{e}_2\cdot A)\wedge(\vec{e}_1\wedge\vec{e}_2) \neq \vec{e}_2\cdot(A \wedge(\vec{e}_1\wedge\vec{e}_2))$---for
example, if $A = \vec{e}_3$ then the
RHS is $\vec{e}_3\wedge\vec{e}_1$ while the
LHS is $0$. There is no doubt that it would be much better
if we could work with a structure like that of a matrix algebra. 
Fortunately  there is a way of doing so.

In what follows it will be useful to look to the exterior product as an
operator acting on $\bigwedge(\mathbb{R}^3)$. 
Let us define the operator $\mathbb{E}(\vec{v}):\bigwedge(\mathbb{R}^3) 
\mapsto \bigwedge(\mathbb{R}^3)$ as 
\begin{equation}
\mathbf{E}(\vec{v})[\Phi] = \vec{v}\wedge\Phi . 
\end{equation}
Let us do the same with the interior product, that is, 
let us define 
the operator $\mathbf{E}^\ast(\vec{v}):\bigwedge(\mathbb{R}^3) 
\mapsto \bigwedge(\mathbb{R}^3) $ as 
\begin{equation}
\mathbf{E}^\ast(\vec{v})[\Phi] = \vec{v}\cdot \Phi . 
\end{equation}
We will also use the compact notation  
\begin{equation}
\mathbf{v} = \mathbf{E}(\vec{v}) 
\end{equation}
and 
\begin{equation}
\mathbf{v}^\ast = \mathbf{E}^\ast(\vec{v}) . 
\end{equation}
In this notation, we have
\begin{equation}
\mathbf{v} = v^i \mathbf{e}_i , \qquad \mathbf{v}^\ast = v^i \mathbf{e}_i^\ast , 
\end{equation}
where we used the summation convention, with 
\begin{equation}
\label{notation.e}
\mathbf{e}_i = \mathbf{E}(\vec{e}_i) , \quad i = 1,2,3, 
\end{equation}
and 
\begin{equation}
\label{notation.e.ast}
 \mathbf{e}^\ast_i = \mathbf{E}^\ast(\vec{e}_i)  , \quad i = 1,2,3. 
\end{equation}
Note that the $\mathbf{e}_i^*$ are operator representations of the dual functionals for the vectors $\vec{e}_i$.

The commutation relations that follow from the skew-symmetry of the
exterior product, from the definition of eq.\eqref{def.int.biv} and 
from eq.\eqref{aux} are 
\begin{gather}
\label{CAR.1.g}
\mathbf{v}\mathbf{u} + \mathbf{u}\mathbf{v} = 0 , \\
\label{CAR.2.g}
\mathbf{v}^\ast \mathbf{u}^\ast + \mathbf{u}^\ast \mathbf{v}^\ast = 0 , \\
\label{CAR.3.g}
\mathbf{v}\mathbf{u}^\ast + \mathbf{u}^\ast \mathbf{v} = 
\vec{v}\cdot \vec{u} . 
\end{gather}
Particular cases are 
\begin{gather} 
\label{CAR.1}
\mathbf{e}_i \mathbf{e}_j + \mathbf{e}_j \mathbf{e}_i = 0 , \\
\label{CAR.2}
\mathbf{e}^\ast_i\mathbf{e}^\ast_j + 
\mathbf{e}^\ast_j\mathbf{e}^\ast_i = 0 , \\
\label{CAR.3}
\mathbf{e}_i \mathbf{e}^\ast_j + \mathbf{e}^\ast_j \mathbf{e}_i = 
\delta_{ij} , 
\end{gather}
for $i, j = 1,2,3$. 

Several readers should have noted at this point the similarity of 
equations \eqref{CAR.1}, \eqref{CAR.2}, \eqref{CAR.3} 
and the anticommutation relations of fermionic creation 
and annihilation operators in quantum theory~\cite{QFT}. In a 
standard notation, if $a_i^\dagger$ and $a_i$ denote the
creation and the annihilation operators of fermionic mode $i$, 
their commutation relations are 
\begin{gather}
\label{ferm1}
a_i^\dagger a_j^\dagger + a_j^\dagger a_i^\dagger = 0 , \\
\label{ferm2}
a_i a_j + a_j a_i = 0 , \\
\label{ferm3}
a_i a_j^\dagger + a_j^\dagger a_i = \delta_{ij} . 
\end{gather}
If the vacuum is denoted by $|0\rangle$, the annihilation 
operators $a_i$ are such that $a_i |0\rangle = 0$. 
The similarity of equations \eqref{ferm1}, \eqref{ferm2}, \eqref{ferm3}, and 
equations \eqref{CAR.1}, \eqref{CAR.2}, \eqref{CAR.3}  is evident, and  
suggests an interesting interpretation for this
formalism. Because of
eq.\eqref{vacuum}, the operators $\mathbf{e}^\ast_i$ plays the
role of annihilation operators and $1$ plays the role of the
vacuum, 
\begin{equation}
\label{annihilation}
\mathbf{e}^\ast_i[1] = 0 , \qquad i = 1,2,3. 
\end{equation}
An arbitrary element of $\bigwedge(\mathbb{R}^3)$ 
can be written as the result of the action of the respective creation 
operator on the vacuum, that is, 
\begin{equation}
\label{act.vac}
(\vec{e}_1)^{\mu_1} \wedge (\vec{e}_2)^{\mu_2} \wedge 
(\vec{e}_3)^{\mu_3}  = 
(\mathbf{e}_1)^{\mu_1}(\mathbf{e}_2)^{\mu_2} 
(\mathbf{e}_3)^{\mu_3}[1] . 
\end{equation}
The idea therefore is to replace vectors by operators according 
to the map 
\begin{equation}
\label{incl}
\imath\big((\vec{e}_1)^{\mu_1} \wedge (\vec{e}_2)^{\mu_2} \wedge 
(\vec{e}_3)^{\mu_3}\big) = (\mathbf{e}_1)^{\mu_1}(\mathbf{e}_2)^{\mu_2} 
(\mathbf{e}_3)^{\mu_3} .
\end{equation}
Let us call eq.\eqref{incl} the {\em natural map}. 
Then we can work with a structure like a matrix algebra, and 
in the end of the calculations, we can get 
the results in terms of vectors using eq.\eqref{annihilation} and eq.\eqref{act.vac}.

We will also need to work 
with the Hodge star operator. 
To write a definition for this operator, let us introduce the following compact notation: 
\begin{gather}
\mathbf{e}_{\mu_1\cdots\mu_k} = \mathbf{e}_{\mu_1}\cdots \mathbf{e}_{\mu_k} , \\
\mathbf{e}^\ast_{\mu_1\cdots\mu_k} = 
\mathbf{e}^\ast_{\mu_1}\cdots \mathbf{e}^\ast_{\mu_k} , \\
\{\mathbf{e}^\ast_i|\mathbf{e}_{\mu_1\cdots\mu_k}\} = 
\mathbf{e}^\ast_i \mathbf{e}_{\mu_1\cdots\mu_k} -
(-1)^k \mathbf{e}_{\mu_1\cdots\mu_k} \mathbf{e}^\ast_i , \\
\{\mathbf{e}^\ast_{\nu_1\cdots\nu_j}|\mathbf{e}_{\mu_1\cdots\mu_k}\} = 
\{\mathbf{e}^\ast_{\nu_1}|\cdots\{\mathbf{e}^\ast_{\nu_j}|
\mathbf{e}_{\mu_1\cdots\mu_k}\}\cdots\} ,
\end{gather}
where $j \leq k$. For example: 
\begin{gather}
\{\mathbf{e}^\ast_i|\mathbf{e}_j\} = \mathbf{e}^\ast_i \mathbf{e}_j + 
\mathbf{e}_j \mathbf{e}^\ast_i , \\
\{\mathbf{e}^\ast_i|\mathbf{e}_j \mathbf{e}_k \} 
= \mathbf{e}^\ast_i \mathbf{e}_j \mathbf{e}_k - \mathbf{e}_j \mathbf{e}_k 
\mathbf{e}^\ast_i , \\
\{\mathbf{e}^\ast_i\mathbf{e}^\ast_j | \mathbf{e}_k \mathbf{e}_l \} = 
\{\mathbf{e}^\ast_i | \{\mathbf{e}^\ast_j | \mathbf{e}_k \mathbf{e}_l \} \} . 
\end{gather}
Let us also denote 
\begin{equation}
\boldsymbol{\Omega} = \mathbf{e}_1 \mathbf{e}_2 \mathbf{e}_3 , 
\end{equation}
and the transformation 
\begin{equation}
\tau(\mathbf{e}_i) = \mathbf{e}_i^\ast, 
\end{equation}
which we generalize as 
\begin{equation}
\tau(\mathbf{e}_{\mu_1\cdots\mu_k}) = \mathbf{e}^\ast_{\mu_1\cdots\mu_k} .
\end{equation}
To define the Hodge star operator acting on $\{\mathbf{e}_i\}$ operators, 
we will look for a generalization of eq.\eqref{hodge.equiv}. 
We define the $\star$ operator  as
\begin{equation}
 \star{\boldsymbol{1}} = \boldsymbol{\Omega} ,\qquad 
 \star{(\mathbf{e}_{\mu_1\cdots\mu_k})} = 
\{\tau\left(\widetilde{\mathbf{e}_{\mu_1\cdots\mu_k}}\right)|\boldsymbol{\Omega}\} . 
\end{equation}
Using these definitions, we have 
\begin{equation*}\begin{aligned}
 \star \boldsymbol{1} &= \mathbf{e}_1 \gp \mathbf{e}_2\gp \mathbf{e}_3 , &
 \star \mathbf{e}_1 &= \mathbf{e}_2 \gp\mathbf{e}_3 , \\
 \star \mathbf{e}_2 &= \mathbf{e}_3\gp\mathbf{e}_1, &
\star \mathbf{e}_3 &= \mathbf{e}_1\gp\mathbf{e}_2 , \\
 \star \mathbf{e}_1\gp\mathbf{e}_2 &= \mathbf{e}_3 , &
 \star \mathbf{e}_3\gp\mathbf{e}_1 &= \mathbf{e}_2 , \\
 \star \mathbf{e}_2\gp\mathbf{e}_3 &= \mathbf{e}_1 , &
 \star \mathbf{e}_1\gp\mathbf{e}_2 \gp\mathbf{e}_3 &= \boldsymbol{1} , 
\end{aligned}\end{equation*}
as expected. We can also define an analogous operation acting on $\{\mathbf{e}^\ast_i\}$ 
operators, as well as on products of $\{\mathbf{e}_i\}$ and 
$\{\mathbf{e}^\ast_i\}$ operators, but we omit these definitions since we 
do not need this operator in this work. 

Now the idea behind our study of the transformations is to 
work with operators instead of vectors. Since the operators
act as linear transformations, the operators are essentially a matrix
algebra. This matrix algebra is the one generated by 
the matrices representing the operators $\mathbf{e}_i$ ($i=1,2,3$) 
and $\mathbf{e}_i^\ast$ ($i=1,2,3$) subject to the
conditions in equations \eqref{CAR.1}, \eqref{CAR.2} and \eqref{CAR.3}. 
Of course we do not need a matrix representation of 
$\mathbf{e}_i$ ($i=1,2,3$) 
and $\mathbf{e}_i^\ast$ ($i=1,2,3$) to study 
the transformations, but if one wants to, it is 
just  a matter of finding the matrix representation
that best serves the application.

\section{Transformation of Points}
\label{sec.5}

Let us consider an arbitrary point $P$, represented
by the paravector $\mathsf{P} = 1 + \vec{p}$. 
Using eq.\eqref{incl}, we write its operator form as 
\begin{equation}
\label{arbitrary.P}
\mathbf{{P}} = \mathbf{1} + \mathbf{p} . 
\end{equation}
We want to study transformations on paravectors of the form
$\mathbf{P} \mapsto  V \mathbf{P} W$, 
and this can be done with the help of the following result. 

\begin{thm}
\label{theorem.0}
Let $\mathbf{P}^\prime = V \mathbf{P} W$, where $\mathbf{P} = \mathbf{1} 
+ \mathbf{p}$ 
is a paravector representing a three dimensional point and $V$ and $W$ 
are elements of the algebra of transformations such that $V\tilde{V}\neq 0$ 
and $W\tilde{W} \neq 0$.  If 
$W = \epsilon \tilde{V}$ ($\epsilon = \pm 1$) and
$V$ is a linear combination of $\mathbf{1}$ and of elements 
that are products of 
terms of the form 
\begin{equation}
\label{form.U}
U = \mathbf{v} \mathbf{u}^\ast ,
\end{equation}
then this transformation is
invertible and $\mathbf{P}^\prime$ is a paravector. 
\end{thm}

\medskip

\begin{proof}
The proof of this theorem 
is left to Appendix~\ref{appendix.A}.
\end{proof}

For $U$ of the form of eq.\eqref{form.U}, its action 
on the vector part of $\mathbf{P}$ is  
\begin{equation}
U \mathbf{p} \tilde{U} = (\vec{p}\cdot\vec{u}) 
(\vec{u}\cdot\vec{v}) \mathbf{v} . 
\end{equation}
A particular case is 
\begin{equation}
U^\prime = \mathbf{v}\mathbf{v}^\ast .
\end{equation}
If we choose the vector $\vec{v}$ such that $|\vec{v}|^2 = 1$, 
then 
\begin{equation}
\label{projection}
U^\prime \mathbf{p}\tilde{U^\prime} = (\vec{p}\cdot\vec{v}) 
\mathbf{v} , 
\end{equation}
which is the projection of $\vec{p}$ in the direction of $\vec{v}$. 
The action of $U$ on the scalar part is $U\tilde{U}$. 
The idea is that different choices of $\mathbf{v}$ and
$\mathbf{u}^\ast$ in eq.\eqref{form.U} and different combinations
of them lead us to the transformations in which we are interested, as we will
see below.

Theorem~\ref{theorem.0}, however, does not exhaust the possible transformations. 
Given an operator $W$, we can calculate $W^n$ and define 
another important class of transformations  
of the form 
\begin{equation}
\mathbf{P} \mapsto \Psi \mathbf{P} \tilde{\Psi} , 
\end{equation}
with 
\begin{equation}
\Psi = {\mbox e}^{t W}  = \sum_{n=0}^\infty \frac{t^n}{n!}W^n .
\end{equation}
Although $W$ can be of the same form as $U$ in eq.\eqref{form.U}, 
we do not need to suppose this---as we will see below---so we 
use a different notation. 
One interesting and important case is when
\begin{equation}
W^2 = a + b W , 
\end{equation}
where $a$ and $b$ are scalars. Then we can write 
\begin{equation}
W^n = c_n + d_n W , 
\end{equation}
where $c_n$ and $d_n$ are also scalar
functions of $a$ and $b$, and consequently 
\begin{equation}
\Psi = C(t) + S(t) W , 
\end{equation} 
where the functions $C(t)$ and $S(t)$ are of the form 
\begin{equation}
C(t) = \sum_{n=0}^\infty \frac{c_n}{n!}t^n , \qquad 
S(t) = \sum_{n=0}^\infty \frac{s_n}{n!} t^n . 
\end{equation}
Then we have
\begin{equation}
\label{trans.exp.form}
\Psi \mathbf{P}\tilde{\Psi} = 
{\mbox e}^{tW}\mathbf{P}{\mbox e}^{t\tilde{W}} = 
C^2(t)\mathbf{P} + C(t)S(t)\left( W \mathbf{P}+\mathbf{P}\tilde{W}\right) + 
S^2(t) W\mathbf{P}\tilde{W} . 
\end{equation}
Note now that we do not need to require $W\mathbf{p}W \neq 0$ 
due to
%eq.\eqref{aux.1a},
eq.\eqref{CAR.1.g},
but only that
$W \mathbf{P}+\mathbf{P}\tilde{W}$ and $W\mathbf{P}\tilde{W}$ 
do not contain $\{\mathbf{e}_i^\ast\}$ operators. Even if 
$W \mathbf{p}+\mathbf{p}\tilde{W} = 0$ and 
$W\mathbf{p}\tilde{W} = 0$, we still have a non-null 
paravector transformation using eq.\eqref{trans.exp.form} 
if $C^2(t) \neq 0$, with translation being an example of this kind of transformation.

There is also another type of transformation that we
can define by exploiting the star operators.  
Given $\mathbf{P}$ as in eq.\eqref{arbitrary.P}, we have
\begin{equation}
\label{star.P}
\star\mathbf{P} = \mathbf{e}_1 \mathbf{e}_2 \mathbf{e}_3 + 
p^1 \mathbf{e}_2\mathbf{e}_3+ p^2 \mathbf{e}_3 \mathbf{e}_1 + 
p^3 \mathbf{e}_1 \mathbf{e}_2 . 
\end{equation}
Let us consider a transformation of the form 
\begin{equation}
\star \mathbf{P} \mapsto \Psi (\star\mathbf{P}) \Phi . 
\end{equation}
Since $\widetilde{\star\mathbf{P}} = -\star\mathbf{P}$, and 
we expect the result of the transformation $\Psi (\star\mathbf{P}) \Phi$ to 
satisfy the same property, we must have  
\begin{equation}
\Phi = \tilde{\Psi} . 
\end{equation}
However, like the discussion in Remark~\ref{remark1}, 
we must also check explicitly that the result of 
this transformation does not contain terms involving $\{\mathbf{e}_i\}$ 
operators. Now we can come back to the space of paravectors 
using $\star$, defining therefore the transformation 
\begin{equation}
\label{trans.star}
\mathbf{P} \mapsto \star (\Psi (\star \mathbf{P})\tilde{\Psi}) . 
\end{equation}

\begin{rem} 
In what follows, we use the following notation for the commutator:
\begin{equation}
[\mathbf{u},\mathbf{v}] = \mathbf{u}\mathbf{v}-\mathbf{v}\mathbf{u} . 
\end{equation}
There are some general relations involving commutators that
are useful. Some of them are 
\begin{gather}
\label{Leibniz}
[\mathbf{u},\mathbf{v}\mathbf{w}] = [\mathbf{u},\mathbf{v}]\mathbf{w} + 
\mathbf{v}[\mathbf{u},\mathbf{w}] , \\
\label{Jacobi}
[\mathbf{u},[\mathbf{v},\mathbf{w}]] + [\mathbf{v},[\mathbf{w},\mathbf{u}]] 
+ [\mathbf{w},[\mathbf{u},\mathbf{v}]] = 0 . 
\end{gather}
Eq.\eqref{Leibniz} is known as the Leibniz rule (since it resembles the 
rule for the derivative of a product) and eq.\eqref{Jacobi} is the 
Jacobi identity.
\end{rem}

\subsection{Reflection}
\label{subs:reflection}

Let us now interpret and identify the generators of some
well-known transformations. 
Consider a point $P$ whose location is specified 
by the vector $\vec{p}$. Given a plane with an unitary 
normal vector $\vec{n}$, the reflection of $\vec{p} = 
\vec{p}_\parallel + \vec{p}_\perp$ 
in this plane generates a vector $\vec{p}{\;}^\prime$ given by 
\begin{equation}
\vec{p}{\;}^\prime =  \vec{p}_\parallel - 
\vec{p}_\perp , 
\end{equation}
where $\vec{p}_\parallel$ is the projection of $\vec{p}$ 
in the mirror plane and $\vec{p}_\perp$ is the component of
$\vec{p}$ in the direction of $\vec{n}$, that is, 
\begin{equation}
\vec{p}_\perp = (\vec{p}\cdot\vec{n}) \vec{n} . 
\end{equation}
The vector $\vec{p}{\;}^\prime$ can be written as
\begin{equation}
\label{reflection.p}
\vec{p}{\;}^\prime = \vec{p} - 2 
(\vec{p}\cdot\vec{n}) \vec{n} .
\end{equation}
Let us see how we can describe this transformation using our algebra.

\begin{thm}
Let $\mathbf{P}$ represent a three dimensional point. 
The point $\mathbf{P}^\prime$ generated by the reflection of
$\mathbf{P}$  
on a plane with a unitary normal vector $\mathbf{n}$ is
given by 
\begin{equation}
\mathbf{P}^\prime = - N \mathbf{P}\tilde{N} , 
\end{equation}
where 
\begin{equation}
N = \mathbf{n}^\ast \mathbf{n} - \mathbf{n}\mathbf{n}^\ast = 
[\mathbf{n}^\ast,\mathbf{n}] . 
\end{equation}
\end{thm}

\medskip

\begin{proof}
First, write 
eq.\eqref{reflection.p} using our operators:
\begin{equation}
  \mathbf{p}^\prime = \mathbf{p} - 2 (\vec{p}\cdot\vec{n})\mathbf{n} 
  = \mathbf{p} + 2 (\vec{p}\cdot\vec{n})\mathbf{n} - 4 (\vec{p}\cdot\vec{n}) 
\mathbf{n} . 
\end{equation}
We have seen that projection can be written in terms of 
eq.\eqref{projection}. So we can write 
\begin{equation}
\mathbf{p}^\prime = 
\mathbf{p} + 2  (\vec{p}\cdot\vec{n})\mathbf{n} -
 4 (\mathbf{n}\mathbf{n}^\ast) 
\mathbf{p} \widetilde{(\mathbf{n}\mathbf{n}^\ast)} . 
\end{equation}
By eq.\eqref{CAR.1.g} and eq.\eqref{CAR.3.g}
\begin{equation}
  \begin{split}
    (\vec{p}\cdot\vec{n})\mathbf{n}&=
    (\mathbf{p}\mathbf{n}^\ast+\mathbf{n}^\ast\mathbf{p})\mathbf{n}\\
    &=\mathbf{p}\mathbf{n}^\ast\mathbf{n}-\mathbf{n}^\ast\mathbf{n}\mathbf{p}\\
    &=\mathbf{p}\mathbf{n}^\ast\mathbf{n}-(1-\mathbf{n}\mathbf{n}^\ast)\mathbf{p}\\
    &=\mathbf{p}\mathbf{n}^\ast\mathbf{n}+\mathbf{n}\mathbf{n}^\ast\mathbf{p}-\mathbf{p}.
  \end{split}
  \end{equation}
Then we have 
\begin{equation}
\begin{split}
\mathbf{p}^\prime & = \mathbf{p} + 2 [\mathbf{p}\mathbf{n}^\ast\mathbf{n}
+ \mathbf{n}\mathbf{n}^\ast\mathbf{p} - \mathbf{p}] - 
4 (\mathbf{n}\mathbf{n}^\ast) \mathbf{p} (\mathbf{n}^\ast\mathbf{n})\\
&= -\mathbf{p} + 2 \mathbf{p}\mathbf{n}^\ast\mathbf{n}
+ 2\mathbf{n}\mathbf{n}^\ast\mathbf{p}  - 
4 (\mathbf{n}\mathbf{n}^\ast) \mathbf{p} (\mathbf{n}^\ast\mathbf{n})\\
& = - (1-2\mathbf{n}\mathbf{n}^\ast) \mathbf{p}(1-2\mathbf{n}^\ast\mathbf{n}) \\
& = - (\mathbf{n}^\ast \mathbf{n} - \mathbf{n}\mathbf{n}^\ast) \mathbf{p} 
(\mathbf{n}\mathbf{n}^\ast - \mathbf{n}^\ast \mathbf{n}) . 
\end{split}
\end{equation}

We also note that  by eq.\eqref{CAR.1.g} and eq.\eqref{CAR.2.g}
\begin{equation}
(\mathbf{n}^\ast \mathbf{n} - \mathbf{n}\mathbf{n}^\ast)
(\mathbf{n}\mathbf{n}^\ast - \mathbf{n}^\ast \mathbf{n}) = 
-\mathbf{n}^\ast\mathbf{n}\mathbf{n}^\ast\mathbf{n} - 
\mathbf{n}\mathbf{n}^\ast\mathbf{n}\mathbf{n}^\ast = 
- \mathbf{n}^\ast\mathbf{n}-\mathbf{n}\mathbf{n}^\ast = - 1. 
\end{equation}
Then the reflected paravector $\mathbf{P}^\prime = \mathbf{1} + 
\mathbf{p}^\prime$ can be written as
\begin{equation}
-\mathbf{P}^\prime = (\mathbf{n}^\ast \mathbf{n} - \mathbf{n}\mathbf{n}^\ast)
(\mathbf{n}\mathbf{n}^\ast - \mathbf{n}^\ast \mathbf{n}) 
+ (\mathbf{n}^\ast \mathbf{n} - \mathbf{n}\mathbf{n}^\ast) \mathbf{p} 
(\mathbf{n}\mathbf{n}^\ast - \mathbf{n}^\ast \mathbf{n}) , 
\end{equation}
that is, 
\begin{equation}
\mathbf{P}^\prime = - N \mathbf{P}\tilde{N} , 
\end{equation}
where 
\begin{equation}
N = \mathbf{n}^\ast \mathbf{n} - \mathbf{n}\mathbf{n}^\ast = 
[\mathbf{n}^\ast,\mathbf{n}] . 
\end{equation}
Note that this is the case where we have $-\tilde{N}$ instead of
$\tilde{N}$ on the RHS of the transformation.
\end{proof}

\subsection{Shear and Non-Uniform Scale Transformations}

The operator $U = \mathbf{u}\mathbf{v}^\ast$ has the following property, which follows directly from eq.\eqref{CAR.3.g},
\begin{equation}
U^2 = \mathbf{u}\mathbf{v}^\ast \mathbf{u}\mathbf{v}^\ast = 
 (\vec{v}\cdot\vec{u})
\mathbf{u}\mathbf{v}^\ast, 
\end{equation}
which generalizes to 
\begin{equation}
U^n = (\vec{v}\cdot\vec{u})^{n-1} \mathbf{u}\mathbf{v}^\ast , \quad 
n = 1,2,3,\ldots 
\end{equation}
This identity suggests that we define a new operator through the exponentiation:
\begin{equation}
{\mbox e}^{t \mathbf{u}\mathbf{v}^\ast} = 1 + 
\sum_{n=1}^\infty \frac{t^n}{n!} (\vec{v}\cdot\vec{u})^{n-1} \mathbf{u}\mathbf{v}^\ast = 1 + H(t) \mathbf{u}\mathbf{v}^\ast ,
\end{equation}
where 
\begin{equation}
\label{eq.def.H}
H(t) = H(t,\vec{v}\cdot\vec{u}) = \begin{cases} 
\displaystyle{\frac{1}{(\vec{v}\cdot\vec{u})} \left(
{\mbox e}^{t(\vec{v}\cdot\vec{u})}-1\right) } , & \quad 
\text{if} \quad (\vec{v}\cdot\vec{u}) \neq 0 , \\
 t ,  & \quad \text{if} \quad 
(\vec{v}\cdot\vec{u}) = 0 .
\end{cases}
\end{equation}
Let us now apply this exponential transformation to a paravector $\mathbf{P}=1+\mathbf{p}$. 
The transformation of the vector part of $\mathbf{P}$ is
\begin{equation}
{\mbox e}^{t \mathbf{u}\mathbf{v}^\ast} \mathbf{p} 
{\mbox e}^{t \mathbf{v}^\ast\mathbf{u}} = 
\mathbf{p} + H(t) \left(\mathbf{u}\mathbf{v}^\ast
\mathbf{p} + \mathbf{p}\mathbf{v}^\ast \mathbf{u} \right) + 
(H(t))^2 \mathbf{u}\mathbf{v}^\ast
\mathbf{p}\mathbf{v}^\ast \mathbf{u} . 
\end{equation}
But by eq.\eqref{CAR.3.g}
\begin{gather}
 \mathbf{u}\mathbf{v}^\ast
\mathbf{p} + \mathbf{p}\mathbf{v}^\ast \mathbf{u} = (\vec{p}\cdot\vec{v}) 
\mathbf{u} + (\vec{u}\cdot\vec{v})\mathbf{p} , \\
\mathbf{u}\mathbf{v}^\ast
\mathbf{p}\mathbf{v}^\ast \mathbf{u}  = (\vec{p}\cdot\vec{v}) 
(\vec{v}\cdot\vec{u})\mathbf{u} , 
\end{gather}
and after some simplifications we obtain 
\begin{equation}
{\mbox e}^{t \mathbf{u}\mathbf{v}^\ast} \mathbf{p} 
{\mbox e}^{t \mathbf{v}^\ast\mathbf{u}} = 
{\mbox e}^{t(\vec{v}\cdot\vec{u})} 
(\mathbf{p} + H(t)
(\vec{p}\cdot\vec{v})\mathbf{u} ) . 
\end{equation}
Moreover, 
\begin{equation}
{\mbox e}^{t \mathbf{u}\mathbf{v}^\ast} \mathbf{1} 
{\mbox e}^{t \mathbf{v}^\ast\mathbf{u}} = 
\mathbf{1}  + H(t)(\vec{v}\cdot\vec{u}) = 
{\mbox e}^{t(\vec{v}\cdot\vec{u})} . 
\end{equation}
The result of this transformation is therefore 
\begin{equation}
{\mbox e}^{t \mathbf{u}\mathbf{v}^\ast} \mathbf{P} 
{\mbox e}^{t \mathbf{v}^\ast\mathbf{u}} = 
{\mbox e}^{t(\vec{v}\cdot\vec{u})}[ \mathbf{1} + 
\mathbf{p} + H(t)
(\vec{p}\cdot\vec{v})\mathbf{u} ] . 
\end{equation}

The factor ${\mbox e}^{t(\vec{v}\cdot\vec{u})}$ 
changes the weight of the point $ \mathbf{1} + 
\mathbf{p} + H(t)
(\vec{p}\cdot\vec{v})\mathbf{u} $ but not 
its location. If we want a transformation that does
not change the weight of the point, we can incorporate 
a factor  
in the transformation that
cancels the factor ${\mbox e}^{t(\vec{v}\cdot\vec{u})}$. 
Obviously this is ${\mbox e}^{-t(\vec{v}\cdot\vec{u})/2}$, 
and the new operator is 
\begin{equation}
{\mbox e}^{-t(\vec{v}\cdot\vec{u})/2} 
{\mbox e}^{t\mathbf{u}\mathbf{v}^\ast} = 
{\mbox e}^{-\frac{t}{2}(\mathbf{u}\mathbf{v}^\ast + \mathbf{v}^\ast 
\mathbf{u}) + t\mathbf{u}\mathbf{v}^\ast} , 
\end{equation}
that is 
\begin{equation}
{\mbox e}^{\frac{t}{2}(\mathbf{u}\mathbf{v}^\ast-\mathbf{v}^\ast\mathbf{u})} = 
{\mbox e}^{t[\mathbf{u},\mathbf{v}^\ast]/2} . 
\end{equation}
It follows that 
\begin{equation}
\label{screw.scale}
{\mbox e}^{t[\mathbf{u},\mathbf{v}^\ast]/2} \mathbf{P}
\widetilde{{\mbox e}^{t[\mathbf{u},\mathbf{v}^\ast]/2}}  = 
\mathbf{1} + 
\mathbf{p} + H(t)
(\vec{p}\cdot\vec{v})\mathbf{u} . 
\end{equation}

Let us now look at two particular cases in detail, 
namely $\vec{v} = \vec{u}$ and 
$\vec{v}\cdot\vec{u} = 0$. Since
the transformation in eq.\eqref{screw.scale} does not 
change the weight of the point, we will focus 
only on the vector part of the paravector.

\begin{thm}
\label{thm:scale}
Let $\mathbf{p}$ and $\mathbf{v}$ represent three dimensional vectors. Then the
transformation 
\begin{equation}
\mathbf{p}^\prime = \Psi \mathbf{p}\tilde{\Psi}
\end{equation}
with
\begin{equation}
\Psi = {\mbox e}^{t[\mathbf{v},\mathbf{v}^\ast]/2}
\end{equation}
is a non-uniform scale transformation of the
component $\mathbf{p}_\parallel$ of $\mathbf{p}$ in the direction of $\mathbf{v}$, that
is, 
\begin{equation}
\Psi \mathbf{p}\tilde{\Psi} = \mathbf{p}_\perp + 
{\mbox e}^{t|\vec{v}|^2} \mathbf{p}_\parallel .
\end{equation}
\end{thm}

\medskip

\begin{proof}
Let us consider $\vec{v} = \vec{u}$ 
in eq.\eqref{screw.scale}. 
In this case 
\begin{equation}
\label{screw.scale.1}
{\mbox e}^{t[\mathbf{v},\mathbf{v}^\ast]/2} \mathbf{p}
\widetilde{{\mbox e}^{t[\mathbf{v},\mathbf{v}^\ast]/2}}  = 
\mathbf{p} + H(t)
(\vec{p}\cdot\vec{v})\mathbf{v} , 
\end{equation}
where from eq.\eqref{eq.def.H}
\begin{equation}
H(t) = \frac{1}{|\vec{v}|^2}({\mbox e}^{t|\vec{v}|^2} - 1) . 
\end{equation}
Let us decompose $\vec{p}$ into the component 
$\vec{p}_\parallel$ in the direction of $\vec{v}$ 
and the component $\vec{p}_\perp$ orthogonal to $\vec{v}$, where 
\begin{equation}
\vec{p}_\parallel = \frac{\vec{p}\cdot\vec{v}}{|\vec{v}|^2} 
\vec{v} . 
\end{equation}
Then eq.\eqref{screw.scale.1} gives
\begin{equation}
{\mbox e}^{t[\mathbf{v},\mathbf{v}^\ast]/2} \mathbf{p}
\widetilde{{\mbox e}^{t[\mathbf{v},\mathbf{v}^\ast]/2}} = 
\mathbf{p}_\perp + 
\mathbf{p}_\parallel + ({\mbox e}^{t|\vec{v}|^2} -1)
\mathbf{p}_\parallel , 
\end{equation}
that is, 
\begin{equation}
{\mbox e}^{t[\mathbf{v},\mathbf{v}^\ast]/2} \mathbf{p}
\widetilde{{\mbox e}^{t[\mathbf{v},\mathbf{v}^\ast]/2}} = 
\mathbf{p}_\perp + 
{\mbox e}^{t|\vec{v}|^2} \mathbf{p}_\parallel ,
\end{equation}
which is a scale transformation in the direction of the
vector $\vec{v}$.
\end{proof}

Since $\Psi$ of Theorem~\ref{thm:scale} leaves $1$ unchanged, $\Psi$ is a non-uniform scale transformation of points.

\medskip

Now decompose $\vec{p}$ as 
\begin{equation}
\label{decomposition.p}
\vec{p} = \vec{p}_\perp + p_u \vec{u} + 
p_v \vec{v} , 
\end{equation}
where 
\begin{equation}
p_u = \frac{\vec{p}\cdot\vec{u}}{|\vec{u}|^2} ,  \qquad 
p_v = \frac{\vec{p}\cdot\vec{v}}{|\vec{v}|^2} ,
\end{equation}
and $\vec{p}_\perp$ orthogonal to $\vec{v}$ and 
$\vec{u}$. Then we have the following result. 

\begin{thm}
  \label{thm:shear}
Let $\mathbf{p}$, $\mathbf{u}$ and $\mathbf{v}$ represent three dimensional vectors 
such that $\vec{u}\cdot\vec{v} = 0$. Then the
transformation 
\begin{equation}
\mathbf{p}^\prime = \Psi \mathbf{p}\tilde{\Psi}
\end{equation}
with
\begin{equation}
\Psi = {\mbox e}^{t[\mathbf{u},\mathbf{v}^\ast]/2}
\end{equation}
is a shear in the plane spanned by $\mathbf{u}$ and $\mathbf{v}$, 
that is, 
\begin{equation}
\Psi \mathbf{p}\tilde{\Psi} = 
\mathbf{p}_\perp + \left(p_u + t |\vec{v}|^2 p_v\right) \mathbf{u} + 
p_v \mathbf{v} . 
\end{equation}
\end{thm}

\medskip

\begin{proof}
  Let us consider in 
eq.\eqref{screw.scale} the situation where $\vec{v}\cdot \vec{u} = 0$. 
In this case we have $H(t) = t$, and from eq.\eqref{screw.scale} we have 
\begin{equation}
{\mbox e}^{t[\mathbf{u},\mathbf{v}^\ast]/2} \mathbf{p}
\widetilde{{\mbox e}^{t[\mathbf{u},\mathbf{v}^\ast]/2}}  = 
\mathbf{p} + t
(\vec{p}\cdot\vec{v})\mathbf{u} . 
\end{equation}
 Then  we obtain 
\begin{equation}
{\mbox e}^{t[\mathbf{v},\mathbf{v}^\ast]/2} \mathbf{p}
\widetilde{{\mbox e}^{t[\mathbf{v},\mathbf{v}^\ast]/2}}  = 
\mathbf{p}_\perp + p_u \mathbf{u} + 
p_v \mathbf{v} + t
p_v |\vec{v}|^2 \mathbf{u} , 
\end{equation}
that is, 
\begin{equation}
{\mbox e}^{t[\mathbf{v},\mathbf{v}^\ast]/2} \mathbf{p}
\widetilde{{\mbox e}^{t[\mathbf{v},\mathbf{v}^\ast]/2}}  = 
\mathbf{p}_\perp + \left(p_u + t |\vec{v}|^2 p_v)\right) \mathbf{u} + 
p_v \mathbf{v} , 
\end{equation}
which is a shear transformation in the plane of the vectors 
$\vec{v}$ and $\vec{u}$.
\end{proof}

Since $\Psi$ of Theorem~\ref{thm:shear} leaves $1$ unchanged, $\Psi$ is a shear transformation of points.

\subsection{Rotations}

We have seen that non-uniform scale and shear transformations are
special cases of the transformation having the operator 
$[\mathbf{u},\mathbf{v}^\ast]$ as its generator. Since 
$[\mathbf{u},\mathbf{v}^\ast] \neq [\mathbf{v},\mathbf{u}^\ast]$, 
we can also define two new operators, namely 
\begin{gather}
\label{eq.R}
\mathcal{R} = [\mathbf{u},\mathbf{v}^\ast] - [\mathbf{v},\mathbf{u}^\ast] , \\
\label{eq.HS}
\mathcal{S} = [\mathbf{u},\mathbf{v}^\ast] + [\mathbf{v},\mathbf{u}^\ast] . 
\end{gather}
The transformations associated with these operators can be summarized 
as follows. 

\begin{thm}
\label{thm.4}
Let $\mathbf{p}$, $\mathbf{u}$ and $\mathbf{v}$ represent three dimensional vectors 
with $\vec{u}\cdot\vec{v} = 0$ and 
$|\vec{u}| = |\vec{v}| = 1$. Then the 
transformation $\Psi\mathbf{p}\tilde{\Psi}$ with 
\begin{equation}
\Psi = {\mbox e}^{\theta \mathcal{R}/2}
\end{equation}
with $\mathcal{R}$ as in eq.\eqref{eq.R} 
is a rotation of the vector $\mathbf{p}$ by an angle $\theta$ in the
plane of $\mathbf{u}$ and $\mathbf{v}$, that is, 
\begin{equation}
\Psi \mathbf{p} \tilde{\Psi} 
 = \mathbf{p}_\perp + 
\mathbf{u}\left[\cos\theta \, p_u + \sin\theta \, p_v\right] + 
\mathbf{v}\left[ \cos\theta \, p_v - \sin\theta \, p_u \right] , 
\end{equation}
where $p_u$ and $p_v$ are defined as in eq.\eqref{decomposition.p}.  
\end{thm}

\medskip

\begin{thm}
\label{thm.5}
Let $\mathbf{p}$, $\mathbf{u}$ and $\mathbf{v}$ represent three dimensional vectors 
with $\vec{u}\cdot\vec{v} = 0$ and 
$|\vec{u}| = |\vec{v}| = 1$. Then the 
transformation $\Psi\mathbf{p}\tilde{\Psi}$ with 
\begin{equation}
\Psi = {\mbox e}^{\theta \mathcal{S}/2}
\end{equation}
with $\mathcal{S}$ as in eq.\eqref{eq.HS} 
is a hyperbolic rotation of the vector $\mathbf{p}$ by an angle $\theta$ in the
plane of $\mathbf{u}$ and $\mathbf{v}$, that is, 
\begin{equation}
\Psi \mathbf{p} \tilde{\Psi} 
= \mathbf{p}_\perp + 
\mathbf{u}[\cosh\theta \, p_u + \sinh\theta \, p_v] + 
\mathbf{v}[\cosh\theta \, p_v + \sinh\theta \, p_u] , 
\end{equation}
where $p_u$ and $p_v$ are defined as in eq.\eqref{decomposition.p}. 
\end{thm}

\medskip

\begin{proof}
The proof of these results can be
done, like the previous cases, by explicit calculations, but are 
longer, so we defer their proofs to Appendix~\ref{appendix.B}.
\end{proof}

\subsection{Translation and Cotranslation}

All the transformations we have studied so far involve 
products of operators like $\mathbf{u}\mathbf{v}^\ast$. 
However, we have seen in the discussion following 
eq.\eqref{trans.exp.form} that we also have 
the possibility of transformations $\Psi$ involving 
a single operator $\mathbf{v}$ such as
\begin{equation}
\label{gen.trans}
\Psi = {\mbox e}^{\mathbf{v}/2} = 1 + \frac{1}{2}\mathbf{v} . 
\end{equation}
We look at two cases of transformations involving this generator. 

\begin{thm}
\label{thm.6}
Let $\mathbf{P}= \boldsymbol{1} +\mathbf{p}$ represent a three dimensional point 
and $\mathbf{v}$ represent a three dimensional vector. 
Then the
transformation 
\begin{equation}
\mathbf{P}^\prime = \Psi \mathbf{P}\tilde{\Psi}
\end{equation}
with $\Psi$ as in eq.\eqref{gen.trans} 
is a translation of the point $\mathbf{P}$ by the
vector $\mathbf{v}$, 
that is, 
\begin{equation}
\Psi \mathbf{P}\tilde{\Psi} = 
\mathbf{P} + \mathbf{v} . 
\end{equation}
\end{thm}

\medskip

\begin{proof}
The action on the vector part of $\mathbf{P}$ is 
\begin{equation}
\left(1 + \frac{1}{2}\mathbf{v}\right)\mathbf{p}\left(1 + \frac{1}{2}\mathbf{v}\right) = 
\mathbf{p} + \frac{1}{2}(\mathbf{v}\mathbf{p}+\mathbf{p}\mathbf{v}) + 
\frac{1}{4}\mathbf{v}\mathbf{p}\mathbf{v} = \mathbf{p} 
\end{equation}
due to eq.\eqref{CAR.1.g}. However, the action on the scalar part of 
$\mathbf{P}$ is 
\begin{equation}
\left(1 + \frac{1}{2}\mathbf{v}\right)\mathbf{1}\left(1 + \frac{1}{2}\mathbf{v}\right) = 
\mathbf{1} + \mathbf{v} . 
\end{equation}
Then we have 
\begin{equation}
{\mbox e}^{\mathbf{v}/2}\mathbf{P} \widetilde{{\mbox e}^{\mathbf{v}/2}} = 
\mathbf{P} + \mathbf{v} , 
\end{equation}
which is the translation of the point $P$ by the vector $\mathbf{v}$.
\end{proof}

\begin{rem}
Note that $\Psi = 1 + \frac{1}{2}\mathbf{v}$ sandwiched on $\mathbf{P}$
 is the identity transformation on the vector part of the paravector of the point $P$;
the contribution to the translation of the point comes from the scalar part 
of the paravector. 
In other words, this transformation of translation does not act on vectors, but 
only on points because of their non-null weight.
Any point can therefore be written as a result of this 
operation acting on the origin of the coordinate system, 
that is,
\begin{equation}
\mathbf{P} = {\mbox e}^{\mathbf{p}/2} \mathbf{O}{\mbox e}^{\mathbf{p}/2} . 
\end{equation}
Since $\mathbf{O} = \mathbf{1}$, we have 
\begin{equation}
\mathbf{P} = \left({\mbox e}^{\mathbf{p}/2}\right)^2 . 
\end{equation}
We have therefore obtained a kind of square root of a point, 
that is, a mathematical object ${\mbox e}^{\mathbf{p}/2}$ whose square gives the mathematical 
object $\mathbf{P}$ used to describe the point.
\end{rem}

\begin{rem}
Similar to~\cite{GS}, our translation of a point $\mathbf{P}$ by a vector $\mathbf{v}$ in Theorem~\ref{thm.6} is a shear in the plane spanned by $1,\mathbf{v}$.
\end{rem}

\begin{thm}
\label{thm.7}
Let $\mathbf{P}$ represent a three dimensional point 
and $\mathbf{v}$ represent a three dimensional vector. 
Then the
transformation 
\begin{equation}
\mathbf{P}^\prime = \star[\Psi(\star \mathbf{P})\tilde{\Psi}] 
\end{equation}
with $\Psi$ as in eq.\eqref{gen.trans} 
satisfies  
\begin{equation}
\label{cotrans}
\star[\Psi(\star \mathbf{P})\tilde{\Psi}] = 
\mathbf{P} + \vec{p}\cdot\vec{v} . 
\end{equation}
\end{thm}

Note that this transformation has the effect 
of giving a weight $\vec{p}\cdot\vec{v}$ 
to the point $\mathbf{P}$. We will call
this transformation \textit{cotranslation} 
in analogy with the definition 
of the codifferential operator, that is, 
the codifferential is the composition of 
duality, differential and duality transformations, 
and the transformation in eq.\eqref{cotrans} is
the composition of duality, translation and 
duality transformations.

\medskip

\begin{proof}
Let us first consider 
the action on the scalar part of the
paravector. We have 
\begin{equation}
{\mbox e}^{\mathbf{v}/2} (\star \mathbf{1}){\mbox e}^{\mathbf{v}/2} = 
\left( 1 + \frac{1}{2}\mathbf{v}\right) \mathbf{e}_1 \mathbf{e}_2
\mathbf{e}_3 \left( 1 + \frac{1}{2}\mathbf{v}\right) = 
\mathbf{e}_1 \mathbf{e}_2
\mathbf{e}_3  
\end{equation} 
because of eq.\eqref{CAR.1.g}. Then 
\begin{equation}
\star[{\mbox e}^{\mathbf{v}/2} (\star \mathbf{1})
 {\mbox e}^{\mathbf{v}/2} ]
= \mathbf{1} ,
\end{equation}
that is, this transformation does not change the weight of the point. However, 
when we consider the vector part of the paravector, we have  
\begin{equation}
{\mbox e}^{\mathbf{v}/2} (\star \mathbf{p}) {\mbox e}^{\mathbf{v}/2} = 
\star\mathbf{p} + \frac{1}{2}(\mathbf{v} (\star \mathbf{p}) + 
(\star \mathbf{p})\mathbf{v}) + \frac{1}{4}\mathbf{v}
(\star\mathbf{p})\mathbf{v} . 
\end{equation}
Eq.\eqref{star.P} gives $\star\mathbf{p}$. 
If we write $\mathbf{v}$ as 
\begin{equation}
\mathbf{v} = v^1 \mathbf{e}_1 + v^2 \mathbf{e}_2 + 
v^3 \mathbf{e}_3
\end{equation}
and use eq.\eqref{CAR.2.g} we find that 
\begin{alignat}{2}
& \mathbf{v} \mathbf{e}_2 \mathbf{e}_3 = 
\mathbf{e}_2 \mathbf{e}_3\mathbf{v} = 
v^1 \mathbf{e}_1 \mathbf{e}_2 \mathbf{e}_3 , & \qquad &
\mathbf{v} \mathbf{e}_3\mathbf{e}_1 = 
\mathbf{e}_3 \mathbf{e}_1 \mathbf{v} = 
v^2 \mathbf{e}_1 \mathbf{e}_2 \mathbf{e}_3 , \\
&\mathbf{v} \mathbf{e}_1 \mathbf{e}_2 = 
\mathbf{e}_1 \mathbf{e}_2 \mathbf{v} = 
v^3 \mathbf{e}_1 \mathbf{e}_2 \mathbf{e}_3 ,& \qquad &
\mathbf{v}\mathbf{e}_1 \mathbf{e}_2 \mathbf{e}_3
 \mathbf{v}= 0 , 
\end{alignat}
and using these expressions in eq.\eqref{star.P} we obtain 
\begin{equation}
{\mbox e}^{\mathbf{v}^/2} (\star \mathbf{p}) {\mbox e}^{\mathbf{v}/2}  = 
\star\mathbf{p} + 
p^1 v^1 \mathbf{e}_1 \mathbf{e}_2 \mathbf{e}_3 + 
p^2 v^2 \mathbf{e}_1 \mathbf{e}_2 \mathbf{e}_3 + 
p^3 v^3 \mathbf{e}_1\mathbf{e}_2\mathbf{e}_3 = 
\star \mathbf{p} + (\vec{p}\cdot\vec{v}) 
\mathbf{e}_1 \mathbf{e}_2 \mathbf{e}_3 . 
\end{equation}
Finally, 
\begin{equation}
\label{weight.vector}
\star[{\mbox e}^{\mathbf{v}/2} (\star \mathbf{p}) 
{\mbox e}^{\mathbf{v}/2} ] = \vec{p}\cdot\vec{v} + \mathbf{p} , 
\end{equation}
which combined with the scalar part gives the result.
\end{proof}

\begin{rem}
Observe that, while the effect of the translation in 
Theorem~\ref{thm.6} comes from the weight of the point, in Theorem~\ref{thm.7}
the effect of the transformation comes from the vector 
part of the point. We can therefore apply this transformation 
to vectors only, as in eq.\eqref{weight.vector}, resulting in a weighted point. 
In other words, 
the transformations in Theorem~\ref{thm.7} transforms a vector into a weighted point with 
weight $\vec{p}\cdot\vec{v}$.
\end{rem}

\begin{rem}
We observe that cotranslation when applied to a point at infinity appears in \cite{Dorst2} under the name perspectivity.  However, cotranslation can be generalized to apply to any object as discussed in Remark~\ref{rem.9}.
\end{rem}

\subsubsection{Perspective Projection Through the Composition of the Translation and the Cotranslation Transformations}

As an application of  the cotranslation transformations, we will see how perspective projection can be described in this formalism. Let us consider two points $P$ and
$E$, described by the paravectors $\mathbf{P} = \mathbf{1} + \mathbf{p}$ and 
$\mathbf{E} = \mathbf{1} + \mathbf{e}$, and a plane with a normal vector $\vec{n}$ and
plane equation 
\begin{equation}
\label{eq.plane}
\vec{x}\cdot\vec{n} = c . 
\end{equation}
To facilitate the discussion, let us introduce the notation 
\begin{equation}
\mathfrak{T}_{\vec{v}}(\mathbf{P}) = {\mbox e}^{\mathbf{v}/2}
\mathbf{P}{\mbox e}^{\mathbf{v}/2} , \qquad 
\mathfrak{W}_{\vec{v}}(\mathbf{P}) = \star[
{\mbox e}^{\mathbf{v}/2}
(\star \mathbf{P}){\mbox e}^{\mathbf{v}/2} ]. 
\end{equation}
Let us start by translating all objects in such a way that the
point $\mathbf{E}$ is moved to the origin. The new points are
\begin{gather}
\mathbf{P}^\prime = \mathfrak{T}_{-\vec{e}}(\mathbf{P}) = 
\mathfrak{T}_{\vec{e}}^{-1}(\mathbf{P}) = 
\mathbf{1} + \mathbf{p}-\mathbf{e}, \\\
\mathbf{E}^\prime = \mathfrak{T}_{-\vec{e}}(\mathbf{E}) = 
\mathfrak{T}_{\vec{e}}^{-1}(\mathbf{E}) = \mathbf{1} . 
\end{gather}
The translated plane equation is 
\begin{equation}
\label{translated.plane}
\vec{y}\cdot\vec{n} = a = c - \vec{n}\cdot\vec{e} , 
\end{equation}
where $\vec{y} = \vec{x}-\vec{e}$. 
Now let us apply the cotranslation transformation by the vector $\vec{n}/a$ 
to the points $\mathbf{P}^\prime$ and $\mathbf{E}^\prime$. 
The results are 
\begin{gather}
\begin{split}
\mathbf{P}^{\prime\prime} = \mathfrak{W}_{\vec{n}/a}(\mathbf{P}^\prime)& = 
\mathbf{1} + \frac{\vec{n}\cdot\vec{p}{\;}^\prime}{a} + \mathbf{p}^\prime \\
& = 
\mathbf{1} +
\frac{\vec{n}\cdot(\vec{p}-\vec{e})}{a} + \mathbf{p}-\mathbf{e} , \end{split} \\ 
\mathbf{E}^{\prime\prime} = \mathfrak{W}_{\vec{n}/a}(\mathbf{E}^\prime) = 
\mathbf{1} . 
\end{gather}
Finally, let us apply the inverse translation to $\mathbf{P}^{\prime\prime}$ and $\mathbf{E}^{\prime\prime}$, to obtain 
\begin{gather}
\label{twt.p}
\begin{split}
\mathbf{P}^{\prime\prime\prime} = \mathfrak{T}_{\vec{e}}(\mathbf{P}^{\prime\prime}) & = 
\left[\mathbf{1}+\frac{\vec{n}\cdot(\vec{p}-\vec{e})}{a}\right](\mathbf{1} + \mathbf{e}) + \mathbf{p}-\mathbf{e} \\ & = 
\mathbf{1} + \frac{\vec{n}\cdot(\vec{p}-\vec{e})}{a} + 
\mathbf{p} + \frac{\vec{n}\cdot(\vec{p}-\vec{e})}{a}\mathbf{e} ,\end{split} \\
\label{twt.e}
\mathbf{E}^{\prime\prime\prime} = \mathfrak{T}_{\vec{e}}(\mathbf{E}^{\prime\prime}) = 
\mathbf{1}+\mathbf{e} = \mathbf{E} . 
\end{gather}
Obviously  eq.\eqref{twt.e} follows from eq.\eqref{twt.p} when we set
$\mathbf{P} = \mathbf{E}$. Let us denote the composition of these
transformations as 
\begin{equation}
\mathfrak{P}_{\vec{e},\vec{n}/a}(\mathbf{P}) = 
(\mathfrak{T}_{\vec{e}}\circ \mathfrak{W}_{\vec{n}/a}
\circ \mathfrak{T}_{\vec{e}}^{-1})(\mathbf{P}) . 
 \end{equation}
From this transformation we have the following, which is similar to~\cite{Goldman3}.

\begin{thm}
Let $\mathbf{n}$ describe the normal to the perspective plane $\mathcal{P}$ 
with equation $\vec{x}\cdot \vec{n}  = c$, the paravector 
$\mathbf{E}$ describe the eye point $E$, and $\mathbf{P}$ describe
an arbitrary point $P$ in three dimensional space. Then 
$\mathbf{P}_0$ given by 
\begin{equation}
\mathbf{P}_0 = 
\mathfrak{P}_{\vec{e},\vec{n}/a}(\mathbf{P}-\mathbf{E})
\end{equation}
where $a = c + \vec{n}\cdot\vec{e}$, 
is a weighted point in the perspective plane located
at the perspective projection of $P$ from the eye point $E$ if
$P$ is located in front of $E$; if $P$ is located behind
$E$, this same location corresponds to 
the weighted point $\bar{\mathbf{P}}_0$. 
\end{thm}

\medskip

\begin{proof}
  The explicit expression 
for $\mathbf{P}_0$ follows from eq.\eqref{twt.p} and eq.\eqref{twt.e}, 
or equivalently, 
\begin{equation}
\begin{split} 
\mathbf{P}_0 & = 
\mathfrak{P}_{\vec{e},\vec{n}/a}(\mathbf{P}-\mathbf{E}) = 
\mathfrak{P}_{\vec{e},\vec{n}/a}(\mathbf{P}) - 
\mathfrak{P}_{\vec{e},\vec{n}/a}(\mathbf{E}) \\
& = \mathbf{1} + \frac{\vec{n}\cdot(\vec{p}-\vec{e})}{a} + 
\mathbf{p} + \frac{\vec{n}\cdot(\vec{p}-\vec{e})}{a}\mathbf{e} - (1+\mathbf{e}) , 
\end{split}
\end{equation}
that is, 
\begin{equation}
\mathbf{P}_0 = \frac{\vec{n}\cdot(\vec{p}-\vec{e})}{a} + 
\mathbf{p} + \left(\frac{\vec{n}\cdot(\vec{p}-\vec{e})}{a}-1\right) \mathbf{e} . 
\end{equation}
This expression for $\mathbf{P}_0$ represents a weighted point whose location is $\mathbf{p}_0 = \frac{\langle \mathbf{P}_0\rangle_1}{|\langle \mathbf{P}_0\rangle_0|}$. 
If $\langle \mathbf{P}_0 \rangle_0 > 0$, then 
\begin{equation}
\mathbf{p}_0 = \frac{a}{\vec{n}\cdot(\vec{p}-\vec{e})}\left[\mathbf{p} + \left(\frac{\vec{n}\cdot(\vec{p}-\vec{e})}{a}-1\right) \mathbf{e}\right] , 
\end{equation} 
which, after some simplifications, can be written as 
\begin{equation}
\mathbf{p}_0 = \left(\frac{c-\vec{n}\cdot\vec{e}}{\vec{n}\cdot(\vec{p}-\vec{e})}\right)\mathbf{p} 
- \left(\frac{c-\vec{n}\cdot\vec{p}}{\vec{n}\cdot(\vec{p}-\vec{e})}\right)\mathbf{e} , 
\end{equation}
which is the expression for the perspective projection of the point $P$ from
the eye point~$E$~\cite{Goldman3}. Clearly $\vec{p}_0\cdot\vec{n} = c$. 
If $\langle \mathbf{P}_0 \rangle_0 < 0$, we obtain the same location 
for the weighted point $\bar{\mathbf{P}}_0$ because 
$\langle \bar{\mathbf{P}}_0 \rangle_1 = 
-\langle \mathbf{P}_0 \rangle_1$. The expression for the weight is 
\begin{equation}
\langle \mathbf{P}_0 \rangle_0 = \frac{a}{\vec{n}\cdot(\vec{p}-\vec{e})} = 
\frac{c-\vec{n}\cdot\vec{e}}{\vec{n}\cdot(\vec{p}-\vec{e})} . 
\end{equation}
To interpret this expression, we remember eq.\eqref{eq.plane} and 
write $c$ in the form $c = \vec{q}\cdot\vec{n}$, where $\vec{q}$ gives  
the location of a 
point $Q$ in the perspective plane. So we have 
\begin{equation}
\langle \mathbf{P}_0 \rangle_0 = 
\frac{\vec{n}\cdot(\vec{q}-\vec{e})}{\vec{n}\cdot(\vec{p}-\vec{e})} = 
\frac{|\overrightarrow{EQ}|\cos{(\angle (\vec{n},\overrightarrow{EQ}))}}
{|\overrightarrow{EP}|\cos{(\angle (\vec{n},\overrightarrow{EP}))}} .
\end{equation}
Points are in front or behind the eye point in relation to the perspective plane 
depending on the sign of $\cos{(\angle (\vec{n},\overrightarrow{EP}))}$, and 
$\langle \mathbf{P}_0 \rangle_0 > 0$ if the sign of $\cos{(\angle (\vec{n},\overrightarrow{EP}))}$ is equal to the sign of 
$\cos{(\angle (\vec{n},\overrightarrow{EQ}))}$, and since $Q$ 
is in the projective plane, $P$ must be
in front of $E$ to have $\langle \mathbf{P}_0 \rangle_0 > 0$.
\end{proof}

\subsubsection{Pseudo-Perspective Projection Through Cotranslation}

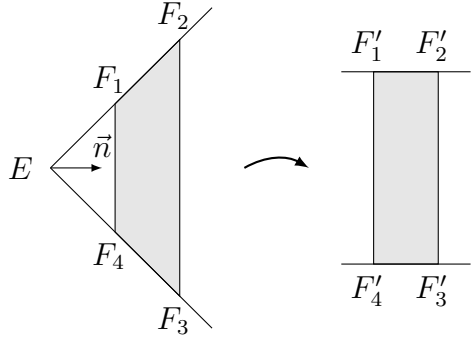
\begin{figure}
\begin{center}
\begin{tikzpicture}[>=latex,xscale = 0.85,yscale=0.85]
\fill [gray!20]
      (1, 1) -- (2,2) -- (2,-2) -- (1,-1) -- cycle;
\draw (1, 1) -- (2,2) -- (2,-2) -- (1,-1) -- cycle;
\draw (0, 0) -- (2.5,2.5);
\draw (0,0) -- (2.5,-2.5);
\node [above] at (0.9,1) {$F_1$};
\node [above] at (1.9,2) {$F_2$};
\node [below] at (1.9,-2) {$F_3$};
\node [below] at (0.9,-1) {$F_4$};
\node [left] at (-0.1,0) {$E$};
\draw [->] (0,0) -- (0.8,0);
\node [above] at (0.8,0) {$\vec{n}$} ;
\draw [thick, ->] (3,0) to [bend left] (4,0);
\fill [gray!20]
      (5, 1.5) -- (6,1.5) -- (6,-1.5) -- (5,-1.5) -- cycle;
      \draw  (5, 1.5) -- (6,1.5) -- (6,-1.5) -- (5,-1.5) -- cycle;
\draw (4.5,1.5) -- (6.5,1.5);
\draw (4.5,-1.5) -- (6.5,-1.5);
\node [above] at (4.9,1.5) {$F_1^\prime$};
\node [above] at (5.9,1.5) {$F_2^\prime$};
\node [below] at (5.9,-1.5) {$F_3^\prime$};
\node [below] at (4.9,-1.5) {$F_4^\prime$};
\end{tikzpicture}
\end{center}
  \caption{Mapping a truncated viewing pyramid to a box.}
  \label{pseudopersp}
\end{figure}

A second application of cotranslation is pseudo-perspective projection.
\textit{Pseudo-perspective projection} is used in computer graphics to map a truncated viewing pyramid (i.e., viewing frustum) to a rectangular box (see Figure~\ref{pseudopersp}); this mapping facilitates z-buffer scan conversion and hidden surface removal by converting the perspective depth test into an orthographic depth test within this box.

Given an eye point $E$ looking in a direction $\vec{n}$, the key observation is that we wish to map the eye point $E$ to a point at infinity, and in particular, we want $E$ to map to $\pm\vec{n}$~\cite{DGM17}.

\begin{thm}
  Let $\mathbf{n}$ be a unit vector, and let $\mathbf{E}=\mathbf{1}-\mathbf{n}$.
  Then
  $\mathfrak{W}_{\vec{n}}(\mathbf{P})$ transforms the eye point $\mathbf{E}$ to the point at infinity in the direction $-\mathbf{n}$ and transforms a viewing frustum to a rectangular box.
\end{thm}

\medskip

\begin{proof}
  Let us start by applying the 
cotranslation operator to $\mathbf{E}$:
\begin{equation}
  \begin{split}
    \mathfrak{W}_{\vec{n}}(\mathbf{E})
    &= \mathfrak{W}_{\vec{n}}(\mathbf{1}-\mathbf{n}) \\
&=\mathbf{1} -\mathbf{n} - \vec{n}\cdot\vec{n} = -\mathbf{n}.
  \end{split}
\end{equation}
Now let $n_{\perp}$ denote a vector perpendicular to $n$.  Consider the four corners of a viewing frustum, as show in Figure~\ref{pseudopersp}. We can represent these corners as 
\begin{alignat*}{2}
& \mathbf{F}_1 = \mathbf{E} + s \mathbf{n} + s \mathbf{n}_\perp , & \qquad & 
\mathbf{F}_2 = \mathbf{E} + t \mathbf{n} + t \mathbf{n}_\perp , \\
& \mathbf{F}_3 = \mathbf{E} + t \mathbf{n} - t \mathbf{n}_\perp , & \qquad & 
\mathbf{F}_4 = \mathbf{E} + s \mathbf{n} - s \mathbf{n}_\perp , 
\end{alignat*}
with $t>s>0$.  
Applying $\mathfrak{W}_{\vec{n}}$ to these points gives the following:
\begin{gather*}
\begin{split}
\mathbf{F}_1^\prime = \mathfrak{W}_{\vec{n}}(\mathbf{F}_1) & = \mathbf{1} - \mathbf{n} + s 
\mathbf{n} + s\mathbf{n}_\perp + (-\vec{n}+ s \vec{n} + s\vec{n}_\perp)\cdot\vec{n} \\
& = s\left(\mathbf{1} + \textstyle{\frac{(s-1)}{s}}\mathbf{n} + 
\mathbf{n}_\perp\right) , 
\end{split} \\
\begin{split}
\mathbf{F}_2^\prime  = \mathfrak{W}_{\vec{n}}(\mathbf{F}_2) & = \mathbf{1} - \mathbf{n} + t 
\mathbf{n} + t\mathbf{n}_\perp + (-\vec{n}+ t \vec{n} + t\vec{n}_\perp)\cdot\vec{n} \\
& = t\left(\mathbf{1} + \textstyle{\frac{(t-1)}{t}}\mathbf{n} + 
\mathbf{n}_\perp \right), 
\end{split} \\
\begin{split}
\mathbf{F}_3^\prime  = \mathfrak{W}_{\vec{n}}(\mathbf{F}_3) & = \mathbf{1} - \mathbf{n} + t 
\mathbf{n} - t\mathbf{n}_\perp + (-\vec{n}+ t \vec{n} - t\vec{n}_\perp)\cdot\vec{n} \\
& = t\left(\mathbf{1} + \textstyle{\frac{(t-1)}{t}}\mathbf{n} - 
\mathbf{n}_\perp\right), 
\end{split} \\
\begin{split}
\mathbf{F}_4^\prime  = \mathfrak{W}_{\vec{n}}(\mathbf{F}_4) & = \mathbf{1} - \mathbf{n} + s 
\mathbf{n} - s\mathbf{n}_\perp + (-\vec{n}+ s \vec{n} - s\vec{n}_\perp)\cdot\vec{n} \\
& = s\left(\mathbf{1} + \textstyle{\frac{(s-1)}{s}}\mathbf{n} - 
\mathbf{n}_\perp \right), 
\end{split} 
\end{gather*}
so that the transformed points are located at the vertices of a rectangle. 
Note that the relative locations of $F_1^\prime$ and $F_2^\prime$ and of 
$F_3^\prime$ and $F_4^\prime$ may be flipped, depending of the original 
locations of the points.
\end{proof}

\begin{rem}
In computer graphics, the normal is usually an axis aligned vector to facilitate hidden surface removal and projection from 3D to 2D.  Further, an arbitrary eye point may be used by translating the arbitrary eye position to $\mathbf{E}=\mathbf{1}-\mathbf{n}$ before performing the pseudo-perspective mapping.
\end{rem}

\section{Transformation of Lines and Planes}
\label{sec.6}

To describe the action of the transformations discussed 
in the previous section on lines and planes, we must first 
incorporate their mathematical description in the formalism. 
We have described a line segment by means of a biparavector 
in eq.\eqref{eq.line}, that is, $\mathcal{L} = 
\mathsf{P}\curlywedge\mathsf{Q}^\dagger = 
\langle \mathsf{P}\wedge\mathsf{Q}^\dagger\rangle_{\{2\}}$. 
The transcription of $\mathcal{L}$ to its operator form 
is straightforward, that is, 
\begin{equation}
\boldsymbol{\mathcal{L}} = \langle \mathbf{P}\mathbf{Q}^\dagger \rangle_{\{2\}} , 
\end{equation}
since the exterior product is already encoded in the definition 
of the operators $\{\mathbf{e}_i\}$. 
Now given the transformations $\mathbf{P}\mapsto \mathbf{P}^\prime = 
U \mathbf{P}V$ and $\mathbf{Q} \mapsto \mathbf{Q}^\prime = 
U\mathbf{Q}V$, with $V = \pm \tilde{U}$, the most natural generalization of its
action to $\boldsymbol{\mathcal{L}}$ is to define 
\begin{equation}
\boldsymbol{\mathcal{L}}^\prime = 
\langle \mathbf{P}^\prime\mathbf{Q}^\prime{}^\dagger \rangle_{\{2\}}.
\end{equation}
Let us work with this expression. 
We have 
\begin{equation}
\boldsymbol{\mathcal{L}}^\prime = \langle U \mathbf{P}V V^\dagger 
\mathbf{Q}^\dagger U^\dagger \rangle_{\{2\}} = 
\langle U \mathbf{P} V\bar{V} \mathbf{Q}^\dagger \bar{U}\rangle_{\{2\}} . 
\end{equation}
Note that 
\begin{equation}
V\bar{V} = \tilde{U}\hat{U} = \widehat{(\bar{U} U)} , 
\end{equation}
and since 
\begin{equation}
\bar{U} U = \overline{\bar{U} U} , 
\end{equation}
we have that $\bar{U} U$ is an element of $\bigwedge^0(\mathbb{R}^3) \oplus
\bigwedge^3(\mathbb{R}^3)$. Let us restrict $U$ so that
$\bar{U} U$ is a scalar. All the transformations studied in the last section
satisfy this property. In fact, they satisfy 
\begin{equation}
\bar{U} U = \varepsilon = \pm 1 , 
\end{equation}
where the case $\varepsilon = -1$ corresponds to the reflection, and all others 
are such that $\varepsilon = 1$. Then we have
\begin{equation}
\boldsymbol{\mathcal{L}}^\prime = \epsilon 
\langle U \mathbf{P} \mathbf{Q}^\dagger \bar{U}\rangle_{\{2\}} . 
\end{equation}
But from the definition of the paravectors $\mathsf{P}$ and 
$\mathsf{Q}^\dagger$ we have that
\begin{equation}
\mathbf{P}\mathbf{Q}^\dagger = \langle \mathbf{P}\mathbf{Q}^\dagger\rangle_{\{2\}} 
- 1, 
\end{equation}
from which we write  
\begin{equation}
U\mathbf{P}\mathbf{Q}^\dagger \bar{U} = U \langle \mathbf{P}\mathbf{Q}^\dagger\rangle_{\{2\}} \bar{U} 
- U\bar{U} = U \boldsymbol{\mathcal{L}} \bar{U} 
- U\bar{U} . 
\end{equation}
Then 
\begin{equation}
\langle U\mathbf{P}\mathbf{Q}^\dagger \bar{U} \rangle_{\{2\}} 
=  \langle U \boldsymbol{\mathcal{L}} \bar{U}  \rangle_{\{2\}} ,
\end{equation}
since $\langle U\bar{U}\rangle_{\{2\}} = 0$. But 
\begin{equation}
\overline{U \boldsymbol{\mathcal{L}} \bar{U} } = 
U \bar{\boldsymbol{\mathcal{L}}} \bar{U} = - U \boldsymbol{\mathcal{L}} \bar{U} ,
\end{equation}
since $\bar{\boldsymbol{\mathcal{L}}} = -\boldsymbol{\mathcal{L}}$ for a 
biparavector $\boldsymbol{\mathcal{L}}$, and then $U \boldsymbol{\mathcal{L}} \bar{U} $
is also a biparavector, 
\begin{equation}
\langle U \boldsymbol{\mathcal{L}} \bar{U}  \rangle_{\{2\}} = 
U \boldsymbol{\mathcal{L}} \bar{U} .
\end{equation}
So we have 
\begin{equation}
\boldsymbol{\mathcal{L}}^\prime = 
\varepsilon U \boldsymbol{\mathcal{L}} \bar{U} ,
\end{equation}
which is the transformation law for biparavectors representing line segments. 

The same line of reasoning can be applied to a triparavector $\mathcal{P}$, 
which is represented in terms of operators as 
\begin{equation}
\boldsymbol{\mathcal{P}} = 
\langle \mathbf{P}\mathbf{Q}^\dagger \mathbf{R}\rangle_{\{3\}}. 
\end{equation}
We obtain that 
\begin{equation}
\boldsymbol{\mathcal{P}}^\prime = \varepsilon 
U \boldsymbol{\mathcal{P}}\tilde{U} 
\end{equation}
is the transformation law for triparavectors representing plane fragments. 
Note that, although the factor $\varepsilon$ 
appears twice from $\bar{U}U$, the quantity $V = \varepsilon \tilde{U}$
appears three times, so in the end we still have a factor $\varepsilon$. 
We also have for a quadriparavector $\boldsymbol{\mathcal{V}}$ that
\begin{equation}
\boldsymbol{\mathcal{V}}^\prime = \varepsilon U \boldsymbol{\mathcal{V}} 
\bar{U} . 
\end{equation}

\begin{rem} 
The fact that
paravectors and triparavectors on one side, and 
biparavectors and quadriparavectors on the other side, 
have different transformation laws is of 
paramount importance. It is the behaviour under 
transformations that enables us to decide whether a
given $k$-vector is part of a $k$-paravector or of
a $(k+1)$-paravector. Let us be more specific: a paravector
is a sum of a scalar and a vector, while a biparavector
is a sum of a vector and a bivector; then, given a 
vector, how does one knows if the vector is part of a paravector
or of a biparavector? The answer comes from the 
transformation properties: if the vector is part of a 
paravector, then the vector transforms as $\mathbf{v} \mapsto 
\varepsilon U \mathbf{v}\tilde{U}$, while if the vector is
part of a biparavector, then the vector transforms as 
$\mathbf{v} \mapsto \varepsilon U \mathbf{v}\bar{U}$. 
However, notice that, with one exception, 
the transformations discussed in last section 
satisfy $\tilde{U} = \bar{U}$. The exception 
is translation, which is generated by 
$\Psi = {\mbox e}^{\mathbf{v}/2}$, and $\bar{\Psi} = {\mbox e}^{-\mathbf{v}/2}$ 
while $\tilde{U} = {\mbox e}^{\mathbf{v}/2}$. 
If $\mathbf{p}$ is part of a paravector, then 
$\mathbf{p}$ transforms as 
\begin{equation}
\mathbf{p} \mapsto {\mbox e}^{\mathbf{v}/2}\mathbf{p} 
{\mbox e}^{\mathbf{v}/2} = \mathbf{p} , 
\end{equation}
that is, $\mathbf{p}$ is not changed. The translation of the points
represented by paravectors comes from the scalar part, which
means that all points of the space with the same weight are
translated equally. On the other hand, if the vector 
is part of a biparavector, $\mathbf{p}$ transforms 
as 
\begin{equation}
\mathbf{p} \mapsto {\mbox e}^{\mathbf{v}/2}\mathbf{p} 
{\mbox e}^{-\mathbf{v}/2} = 
{\mbox e}^{\mathbf{v}}\mathbf{p} = 
\mathbf{p} + \mathbf{v}\mathbf{p} .
\end{equation}
which represents a line segment with direction 
defined by the vector $\vec{p}$, moment 
$\vec{v}\wedge\vec{p}$ about the origin, and  
$\vec{d} = (\vec{v}\wedge\vec{p})\cdot\vec{p}/
|\vec{p}|^2 = \vec{v}-(\vec{v}\cdot\vec{p})\vec{p}/
|\vec{p}|^2$. In other words: 
a line passing through the origin is translated to a parallel
line passing through the point $\mathbf{V} = 1 + \vec{v}$.
\end{rem}

\medskip

\begin{rem}
\label{rem.9}
Now that we know how to apply the translation operator $\mathfrak{T}_{\vec{v}}$ 
to an arbitrary $k$-paravector, we can extend the definition 
of cotranslation to an arbitrary $k$-paravector. 
Let us write the operator form of a $k$-paravector as 
\begin{equation}
\boldsymbol{\mathcal{P}}_{\{k\}} = \boldsymbol{A}_{k-1} + 
\boldsymbol{A}_k . 
\end{equation} 
Then we can show that 
\begin{equation}
\boldsymbol{\mathcal{P}}^\prime_{\{k\}} = 
\star \mathfrak{T}_{\vec{v}}(\star \boldsymbol{\mathcal{P}}_{\{k\}}) = 
\boldsymbol{A}_{k-1} + \boldsymbol{A}_k\cdot \mathbf{v} + \boldsymbol{A}_k = 
\boldsymbol{\mathcal{P}}_{\{k\}} + \boldsymbol{A}_k\cdot \mathbf{v}, 
\end{equation}
which is the generalization of eq.\eqref{cotrans} for an
arbitrary $k$-paravector. 
\end{rem}

\section{Conclusions}
\label{sec.7}

We have provided an intrinsic approach to the geometry of a three dimensional 
Euclidean space and its geometric transformations  
based on an algebra constructed from
elements of the three dimensional space. 
The concept of a paravector was introduced as an algebraic representative
of a point, and a paravector contains information about the location, weight and
orientation of this point. We have introduced a product of paravectors 
giving a biparavector, and when this product involves paravectors 
representing points with opposite orientations, the biparavector 
represents the line segment joining these points, in such a way 
that this biparavector resembles the Pl\"ucker representation of a line. 
The same construction can be applied to the product of three
paravectors, resulting in a triparavector representing a 
plane fragment. Although we have discussed only three
dimensional Euclidean space, this formalism is not restricted
to three dimensions or to the Euclidean case, and can be 
easily generalized. 

We have studied geometric transformations on this
three dimensional space by means of an algebra of 
transformations. We have shown that this formalism describes  
reflection, rotations (circular 
and hyperbolic), translation, shear and non-uniform 
scale transformations in an unified way. 
Using the concept of Hodge duality, 
we have also defined a new operation called cotranslation, and showed 
that the operation of perspective projection can be 
written as a composition of the translation and cotranslation 
operations. 

We have also discussed the subtle difference 
in the transformations of points, lines and planes 
for the case of translations. This difference makes it 
possible to distinguish when a $k$-vector is part
of a $k$-paravector or a $(k+1)$-paravector. 

Many readers must have noticed a relationship between the
algebra of transformations in eq.\eqref{CAR.1}, eq.\eqref{CAR.2} and
eq.\eqref{CAR.3} and the algebraic relations defining Clifford algebras.
In fact, we can define Clifford algebras from those expressions \cite{VazRocha}, 
in particular the Clifford algebras of quadratic spaces with
signature $(3,0)$, $(0,3)$ and $(3,3)$. However, a 
discussion of these Clifford algebras and 
their potentialities in dealing with the various geometric 
transformations is beyond the scope of this work, 
and will be done elsewhere.  
Notwithstanding, we believe some comments about the use of Clifford
algebra in computer graphics are welcome.

Others have used Clifford algebras to create models of affine spaces
for use in computer graphics including perspective projections.
Gunn's~\cite{Gunn} $P(R^*(3,0,1))$ model has some features similar to our model, but lacks shears and non-uniform scaling.
Our model is similar to the $R(4,4)$ model of 
Du et al.~\cite{DGM17,GS}: their $R(4,4)$ model has all the objects and transformations described in this paper, and in addition has a representation for quadric surfaces.  One significant difference is that their model uses an extra dimension in the vector space, resulting in a Clifford algebra four times the size of our model. 
Dorst~\cite{Dorst2} develops a model for the study of oriented projective transformations of lines; our representation of lines is similar to that of Dorst
but ours is a model of affine space rather than focused on lines, and
our model can be generalized to arbitrary dimensions.  A deeper
analysis of both approaches should be done after we formulate a 
version of our work using the same algebra used in~\cite{Dorst2}. 
We also observe that the use of a non-homogeneous combination 
of algebraic elements like in our definition 
of $k$-paravectors was used by Selig~\cite{Selig} in the description 
of some configurations of points and lines and of lines and planes 
in robotics, which he called flags.

The advantages of our approach is that our model contains points,
line segments, and plane sectors in a natural way, and it includes all
affine and projective transformations.  The derivations of shear and
non-uniform scaling are easy, as are the derivations of translation
and cotranslation.  Further, our model includes perspective and
pseudo-perspective, using cotranslation.  The disadvantages of our
approach include that the derivation of reflection is more complicated
than in most competing approaches, and the derivations of rotation and
hyperbolic rotation are much more complicated.

While we used gaigen~\cite{DF} to implement the geometry and verify
the formulas in this paper,
we have not yet considered how much the formalism 
developed in this work can be useful from a practical 
point of view, especially when we think of its applications 
in computer graphics. As we hope to have made clear 
in the introduction, our interest in this work is 
essentially theoretical. However, the continuation of 
this work has to go through this discussion, and we 
intend to do so in a timely manner. We believe that 
the relationship of the formalism presented here with the
formalism of Clifford algebras may be the path that 
could lead to an efficient procedure for practical 
applications of the concepts and results discussed here.

\bigskip

\noindent \textbf{\sffamily Acknowledgements:} JV gratefully 
acknowledges the support of a research grant from 
FAPESP -- process 2016/21370-9.  SM is grateful for the support of the
Natural Sciences and Engineering Research Council of Canada. 
JV is grateful to the University of Waterloo for the hospitality 
during his stay as visiting professor. 

\appendix 

\section{Proof of Theorem~\ref{theorem.0}}
\label{appendix.A}

We want to study transformations on paravectors of the form
$\mathbf{P} \mapsto  V \mathbf{P} W$ and 
we expect the result of this transformation to be another 
paravector, and paravectors are the only elements in three
dimensions that satisfy $\phi = \tilde{\phi}$---but see Remark~\ref{remark1} below. Then, using the property $\widetilde{AB} = \tilde{B}\tilde{A}$, 
the transformation has to satisfy 
\begin{equation}
V \mathbf{P} W = \tilde{W} \mathbf{P} \tilde{V} . 
\end{equation}
Moreover, we also expect the transformation to be 
invertible. Then this transformation also has to satisfy 
\begin{equation}
V \mathbf{P} W = \tilde{W} \mathbf{P} \tilde{V} . 
\end{equation}
Moreover, we also expect this transformation to be 
invertible. Then this transformation also has to satisfy 
\begin{equation}
  \label{eq:satisfy}
\mathbf{P} =  V^{-1} (V \mathbf{P} W) W^{-1} = 
V^{-1} \tilde{W} \mathbf{P} \tilde{V} W^{-1} , 
\end{equation}
and since \eqref{eq:satisfy} has to be valid for all 
paravectors, we conclude that
\begin{equation}
W = \varepsilon\tilde{V} , \quad \varepsilon = \pm 1 . 
\end{equation}
Although we do not disregard the possibility of the minus sign in this transformation---there is indeed an important
example of this case, as we see in Section~\ref{subs:reflection}---in the 
following discussion we will use only the plus sign.

\begin{rem}
\label{remark1}
The algebra of 
transformations has generators $\{\mathbf{e}_i\}$ ($i=1,2,3$) 
and $\{\mathbf{e}^\ast_i\}$ ($i=1,2,3$), so the dimension of this algebra 
is $2^3 \cdot 2^3 = 2^6$, instead of $2^3$ as for $\bigwedge(\mathbb{R}^3)$, 
and an arbitrary element of this algebra can be written 
using equations \eqref{CAR.1}, \eqref{CAR.2} and \eqref{CAR.3} in 
the form 
\begin{equation}
\label{gen.elem.trans}
(\mathbf{e}_1)^{\mu_1}(\mathbf{e}_2)^{\mu_2}
(\mathbf{e}_3)^{\mu_3}(\mathbf{e}_1^\ast)^{\nu_1} 
(\mathbf{e}_2^\ast)^{\nu_2}(\mathbf{e}_3^\ast)^{\nu_3} . 
\end{equation} 
So inside the algebra of transformations, 
the condition $\phi = \tilde{\phi}$ ---see eq.\eqref{eq.tilde}---
implies that 
$\phi$ is a combination of elements of the form given 
by eq.\eqref{gen.elem.trans} with 
$|\mu|+|\nu| = 0,1,4,5$, 
where $|\mu| = \mu_1 + \mu_2 + \mu_3$ 
and $|\nu| = \nu_1 + \nu_2 + \nu_3$. 
These results mean  that we also have
to impose the condition $|\nu| = 0$  
to guarantee that $|\mu| = 0$ (scalar) and 
$|\mu| = 1$ (vector), that is, for the result 
to be a paravector. Thus, the condition 
$|\nu| = 0$---which means that there is no term 
involving the $\{\mathbf{e}_i^\ast\}$ operators in the final 
result---is a restriction on the operators. 
\end{rem}

Let us look for an expression for $V$ as a combination 
of some basic terms, that we will denote by $U$, 
from the algebra of transformation. 
The transformation $\mathbf{P} \mapsto U \mathbf{P}\tilde{U}$ changes
the mass and location of the point according to 
\begin{equation}
\mathbf{1} \mapsto U\tilde{U} , \qquad 
\mathbf{p} \mapsto U \mathbf{p} \tilde{U} . 
\end{equation}
We expect that $U \mathbf{p} \tilde{U} \neq 0$ because 
otherwise all points are mapped to the same point $U\tilde{U}$ and
the transformation cannot be invertible. So let us focus our
attention on the transformation $\mathbf{p} \mapsto U \mathbf{p} \tilde{U}$ 
looking for $U$ such that $ U \mathbf{p} \tilde{U} \neq 0$ 
and with no terms in the final result 
involving the $\{\mathbf{e}_i^\ast\}$ operators. 

From equations \eqref{CAR.1.g}, \eqref{CAR.2.g} and \eqref{CAR.3.g}, we have
that 
\begin{align}
\label{aux.1a}
\mathbf{v}\mathbf{u}\mathbf{v} &= 0 , \\
\label{aux.1b}
\mathbf{v}\mathbf{u}^\ast \mathbf{v} &= (\vec{v}\cdot\vec{u}) \mathbf{v} , \\
\label{aux.1c} 
\mathbf{v}^\ast \mathbf{u}\mathbf{v}^\ast &= 
(\vec{v}\cdot\vec{u})\mathbf{v}^\ast , \\
\label{aux.1d}
\mathbf{v}^\ast \mathbf{u}^\ast \mathbf{v}^\ast &= 0 . 
\end{align}
An arbitrary $U$ (with $U \neq \mathbf{1}$) is a sum 
of products of elements
$\mathbf{v}_1$, $\mathbf{v}_2, \ldots$ and 
$\mathbf{u}_1^\ast$, $\mathbf{u}_2^\ast, \ldots$ 
Then $U$ must be of the form $U = U_1 \mathbf{u}_1^\ast$, 
because if $U$ is of the form $U = U^\prime_1 \mathbf{v}_1$ 
then $U \mathbf{p}\tilde{U} = 0$ because of eq.\eqref{aux.1a}. 
So from eq.\eqref{aux.1c}, we have 
\begin{equation}
U \mathbf{p} \tilde{U} = (\vec{p}\cdot\vec{u}_1) 
U_1 \mathbf{u}_1^\ast \tilde{U}_1 . 
\end{equation}
Now $U_1 \neq \mathbf{1}$ because $U_1 \mathbf{u}_1^\ast \tilde{U}_1$
 involves $\{\mathbf{e}_i^\ast\}$ operators. 
Moreover, $U_1$ cannot be of the form $U_2^\prime \mathbf{u}_2^\ast$ 
because in this case $U_1 \mathbf{u}_1^\ast \tilde{U}_1  = 0$ 
due to eq.\eqref{aux.1d}. Then $U_1$ must be of the form 
$U_1 = U_2 \mathbf{v}_1$, and then, using eq.\eqref{aux.1b},  
\begin{equation}
U\mathbf{p}\tilde{U} = (\vec{p}\cdot\vec{u}_1) 
(\vec{u}_1\cdot\vec{v}_1) U_2 \mathbf{v}_1 \tilde{U}_2 , 
\end{equation}
where
\begin{equation}
U = U_2 \mathbf{v}_1\mathbf{u}^\ast_1 . 
\end{equation}
Now we can have ${U}_2 = \mathbf{1}$, and 
if  ${U}_2 \neq \mathbf{1}$, then 
repeating the preceding discussion, we conclude that
\begin{equation}
U = U_4 \mathbf{v}_2 \mathbf{u}_2^\ast \mathbf{v}_1 \mathbf{u}_1^\ast , 
\end{equation}
and so on. 

We conclude therefore that $V$ is a 
linear combination of $\mathbf{1}$ and of products of terms of the form 
\begin{equation}
U = \mathbf{v} \mathbf{u}^\ast ,
\end{equation}
with 
\begin{equation}
U \mathbf{p} \tilde{U} = (\vec{p}\cdot\vec{u}) 
(\vec{u}\cdot\vec{v}) \mathbf{v} . 
\end{equation}

We also note that an arbitrary term must have the same number of operators
$\mathbf{v}_i$ and $\mathbf{u}_i^\ast$ irrespective of the order, because of
equations \eqref{CAR.1.g}, \eqref{CAR.2.g} and \eqref{CAR.3.g}---for 
example: the term $\mathbf{v}_1\mathbf{v}_2\mathbf{u}_1^\ast \mathbf{u}_2^\ast$ 
can also be written as the sum 
$(\vec{v}_2\cdot\vec{u}_1) 
\mathbf{v}_1\mathbf{u}_2^\ast - \mathbf{v}_1\mathbf{u}_1^\ast 
\mathbf{v}_2\mathbf{u}_2^\ast$, which is of the form of eq.\eqref{form.U}.

\section{Proof of Theorem~\ref{thm.4} and Theorem~\ref{thm.5}}
\label{appendix.B}

To see the result of the transformation ${\mbox e}^{tX/2}
\mathbf{P} {\mbox e}^{t\tilde{X}/2}$ for $X = \{\mathcal{R},\mathcal{S}\}$, 
we need to calculate $X^n$ for $x=2,3,4,\ldots$ It is not difficult to see 
that 
\begin{equation}
\mathcal{R}^2 = 4[(\vec{u}\cdot\vec{v})(\mathbf{u}
\mathbf{v}^\ast + \mathbf{v}\mathbf{u}^\ast) - |\vec{v}|^2 \mathbf{u}
\mathbf{u}^\ast - |\vec{u}|^2\mathbf{v}\mathbf{v}^\ast] - 
8 \mathbf{u}\mathbf{v}\mathbf{u}^\ast\mathbf{v}^\ast 
\end{equation}
and 
\begin{equation}
\mathcal{S}^2 = 4(\vec{u}\cdot \vec{v})^2 - \mathcal{R}^2 . 
\end{equation}
As we see, expressions for $\mathcal{R}^2$ and $\mathcal{S}^2$ are
not simple, so we expect not to have a simple expression for 
${\mbox e}^{tX/2}$ for $X = \{\mathcal{R},\mathcal{S}\}$. 
A strategy to overcome this difficulty is to write (if possible) 
the operator $X$ as 
\begin{equation}
X = X_1 + X_2 
\end{equation}
with 
\begin{equation}
[X_1,X_2] = 0 . 
\end{equation}
If this condition is satisfied, then
\begin{equation}
{\mbox e}^{tX/2} = {\mbox e}^{tX_1/2}{\mbox e}^{tX_2/2} , 
\end{equation}
and if the expressions for $(X_1)^n$ and $(X_2)^n$ are 
not as complicated as the one for $X^2$, we can
obtain expressions for ${\mbox e}^{tX_1/2}$ and ${\mbox e}^{tX_2/2}$ 
that are not so difficult to handle, and then 
study the effect of ${\mbox e}^{tX}$ through this 
decomposition. This is what we do in the following proof of Theorem~\ref{thm.4}.

\begin{proof}
If 
we add and subtract the terms $\frac{1}{2}[\mathbf{u},\mathbf{v}]$ and
$\frac{1}{2}[\mathbf{u}^\ast,\mathbf{v}^\ast]$ to $\mathcal{R}$ as in 
eq.\eqref{eq.R}, we can write
\begin{equation}
\mathcal{R} = \mathcal{R}_1 - \mathcal{R}_2 , 
\end{equation}
where 
\begin{equation}
\label{eq.R1.R2}
\mathcal{R}_1 = \frac{1}{2}[\mathbf{u}+\mathbf{u}^\ast, 
\mathbf{v}+\mathbf{v}^\ast] , \qquad 
\mathcal{R}_2 = \frac{1}{2}[\mathbf{u}-\mathbf{u}^\ast,
\mathbf{v}-\mathbf{v}^\ast] . 
\end{equation}
Now let us calculate $[\mathcal{R}_1,\mathcal{R}_2]$. After using
the property $[X,Y] + [Y,X] = 0$ to make some simplifications, we 
obtain that 
\begin{equation}
[\mathcal{R}_1,\mathcal{R}_2] = 2\big[ [\mathbf{u},\mathbf{v}^\ast] + 
[\mathbf{u}^\ast,\mathbf{v}] , [\mathbf{u},\mathbf{v}] + [\mathbf{u}^\ast,
\mathbf{v}^\ast] \big] . 
\end{equation}
If we write the inside commutators like $[\mathbf{u},\mathbf{v}^\ast] = 
\mathbf{u}\mathbf{v}^\ast - \mathbf{v}^\ast\mathbf{u}$ and use 
\begin{xalignat}{2}
& [\mathbf{u}\mathbf{v}^\ast,\mathbf{u}\mathbf{v}] = 
[\mathbf{v}\mathbf{u}^\ast,\mathbf{u}\mathbf{v}] = (\vec{v}\cdot\vec{u}) 
\mathbf{u}\mathbf{v} , &  & 
 [\mathbf{u}\mathbf{v}^\ast,\mathbf{u}^\ast\mathbf{v}^\ast] = 
[\mathbf{v}\mathbf{u}^\ast,\mathbf{u}^\ast\mathbf{v}^\ast] = -(\vec{v}\cdot\vec{u}) 
\mathbf{u}^\ast\mathbf{v}^\ast , \\
& [\mathbf{u}^\ast\mathbf{v},\mathbf{u}\mathbf{v}] = 
[\mathbf{v}^\ast\mathbf{u},\mathbf{u}\mathbf{v}] = -(\vec{v}\cdot\vec{u}) 
\mathbf{u}\mathbf{v} ,&  & 
[\mathbf{u}^\ast\mathbf{v},\mathbf{u}^\ast\mathbf{v}^\ast] = 
[\mathbf{v}^\ast\mathbf{u},\mathbf{u}^\ast\mathbf{v}^\ast] = (\vec{v}\cdot\vec{u}) 
\mathbf{u}^\ast\mathbf{v}^\ast ,
\end{xalignat}
it follows that 
\begin{equation}
[\mathcal{R}_1,\mathcal{R}_2] = 0 . 
\end{equation}
Note that this result does not depend on the assumption that
$(\vec{v}\cdot\vec{u}) = 0$, which we will use below. 

Next we calculate $(\mathcal{R}_1)^2$ and $(\mathcal{R}_2)^2$. 
The calculation is straightforward but long and tedious, so we leave
the details for Appendix~\ref{appendix.C}.  The result is 
\begin{equation}
(\mathcal{R}_1)^2 = (\mathcal{R}_2)^2 = -|\vec{u}\wedge\vec{v}|^2 
< 0 ,  
\end{equation}
where we are supposing only that the vectors $\vec{u}$ and 
$\vec{v}$ are such that $\vec{u}\wedge\vec{v}\neq0$, 
that is, $\vec{u}$ and $\vec{v}$ are linearly independent. 
Then we conveniently define 
\begin{equation}
\mathcal{I}_{1} = \frac{\mathcal{R}_{1}}{|\vec{u}\wedge\vec{v}|} , 
\qquad \mathcal{I}_{2} = \frac{\mathcal{R}_{2}}{|\vec{u}\wedge\vec{v}|} ,
\end{equation}
in such a way that
\begin{equation}
(\mathcal{I}_{1})^2 = (\mathcal{I}_{2})^2 =  -\mathbf{1} . 
\end{equation}
As a consequence, 
\begin{gather}
\label{euler.1}
{\mbox e}^{\theta \mathcal{I}_{1}/2} = 
\cos{(\theta/2)} + \mathcal{I}_{1} \sin{(\theta/2)} , \\
\label{euler.2}
{\mbox e}^{\theta \mathcal{I}_{2}/2} = 
\cos{(\theta/2)} + \mathcal{I}_{2} \sin{(\theta/2)}  ,
\end{gather}
and we can write 
\begin{equation}
{\mbox e}^{t\mathcal{R}/2} = {\mbox e}^{\theta \mathcal{I}_1/2} 
{\mbox e}^{-\theta\mathcal{I}_2/2} = 
{\mbox e}^{-\theta \mathcal{I}_2/2} 
{\mbox e}^{\theta\mathcal{I}_1/2} ,
\end{equation}
where we identified $t |\vec{u}\wedge\vec{v}| = \theta$, 
which is useful when we have arbitrary vectors (such 
that $\vec{u}\wedge\vec{v}\neq 0$). 
However, when $\vec{u}\cdot\vec{v} = 0$ and $|\vec{v}| = 
|\vec{u}| = 1$, we have $|\vec{u}\wedge\vec{v}| = 1$, 
and then $t = \theta$. 

Since this transformation does not change the weight of a
paravector, let us calculate its effect on an arbitrary vector. 
Let us first calculate 
\begin{equation}
\mathbf{p}^\prime = {\mbox e}^{\theta\mathcal{I}_1/2} \mathbf{p} 
\widetilde{{\mbox e}^{\theta\mathcal{I}_1/2}} = 
{\mbox e}^{\theta\mathcal{I}_1/2} \mathbf{p} 
{\mbox e}^{-\theta\mathcal{I}_1/2} . 
\end{equation}
Using eq.\eqref{euler.1} we have
\begin{equation}
\label{aux.p.prime}
\mathbf{p}^\prime = \cos^2\frac{\theta}{2} + 
\cos\frac{\theta}{2}\sin\frac{\theta}{2}[\mathcal{I}_1,\mathbf{p}] 
-\sin^2\frac{\theta}{2}\mathcal{I}_1\mathbf{p}\mathcal{I}_1 . 
\end{equation}
To calculate $[\mathcal{I}_1,\mathbf{p}]$ 
it is convenient to use some results involving commutators, 
which can be proved using the
commutation relations in equations \eqref{CAR.1.g}, \eqref{CAR.2.g} and
\eqref{CAR.3.g}. Some of these formulas are
\begin{xalignat}{2}
&\left[ [\mathbf{u},\mathbf{v}^\ast], 
\mathbf{p}\right] = (\vec{p}\cdot\vec{v})\mathbf{u} , &  &
\left[ [\mathbf{u},\mathbf{v}^\ast],
\mathbf{p}^\ast \right] = -(\vec{p}\cdot\vec{u})\mathbf{v}^\ast , \\
& \left[ [\mathbf{u},\mathbf{v}], 
\mathbf{p}\right] = 0 , & &
\left[ [\mathbf{u},\mathbf{v}], 
\mathbf{p}^\ast\right] = (\vec{p}\cdot\vec{v})\mathbf{u} - 
(\vec{p}\cdot\vec{u})\mathbf{v} , \\
& \left[ [\mathbf{u}^\ast,\mathbf{v}^\ast], 
\mathbf{p}\right] = (\vec{p}\cdot\vec{v})\mathbf{u}^\ast - 
(\vec{p}\cdot\vec{u})\mathbf{v}^\ast , & & 
\left[ [\mathbf{u}^\ast,\mathbf{v}^\ast], 
\mathbf{p}^\ast\right] = 0.
\end{xalignat}
To calculate $\mathcal{I}_1\mathbf{p}\mathcal{I}_1 $ 
it is convenient to use the following: 
\begin{gather}
(\mathbf{u}+\mathbf{u}^\ast)[\mathbf{u}+\mathbf{u}^\ast,
\mathbf{v}+\mathbf{v}^\ast] = 2 (\mathbf{v} + \mathbf{v}^\ast)  , \\
(\mathbf{v}+\mathbf{v}^\ast)[\mathbf{u}+\mathbf{u}^\ast,
\mathbf{v}+\mathbf{v}^\ast] = -2 (\mathbf{u} + \mathbf{u}^\ast)  , 
\end{gather}
where we used the assumptions that $\vec{v}\cdot \vec{u} = 0$ 
and $|\vec{v}|^2 = |\vec{u}|^2 = 1$. 
With these expressions we obtain 
\begin{equation}
[\mathcal{I}_1,\mathbf{p}] = \frac{1}{|\vec{u}\wedge\vec{v}|} 
\big( (\vec{p}\cdot\vec{v})(\mathbf{u}+\mathbf{u}^\ast) - 
(\vec{p}\cdot\vec{u})(\mathbf{v}+\mathbf{v}^\ast) \big) 
\end{equation}
and
\begin{equation}
\mathcal{I}_1\mathbf{p}\mathcal{I}_1  = 
-\mathbf{p} + \frac{1}{|\vec{u}\wedge\vec{v}|^2} 
\big[ (\vec{p}\cdot\vec{v}) (\mathbf{v} + 
\mathbf{v}^\ast)  
+ (\vec{p}\cdot\vec{u}) (\mathbf{u} + 
\mathbf{u}^\ast)\big] . 
\end{equation}
and using these results in eq.\eqref{aux.p.prime} we conclude that
\begin{equation}
\mathbf{p}^\prime = \mathbf{p} + C_u (\mathbf{u}+\mathbf{u}^\ast) + 
C_v (\mathbf{v}+\mathbf{v}^\ast) 
\end{equation}
where 
\begin{gather}
C_u = \cos\frac{\theta}{2}\sin\frac{\theta}{2} 
\frac{(\vec{p}\cdot\vec{v})}{|\vec{u}\wedge\vec{v}|} - 
\sin^2\frac{\theta}{2} \frac{ (\vec{p}\cdot\vec{u})}{|\vec{u}\wedge\vec{v}|^2} , \\
C_v = -\cos\frac{\theta}{2}\sin\frac{\theta}{2} 
\frac{(\vec{p}\cdot\vec{u})}{|\vec{u}\wedge\vec{v}|} - 
\sin^2\frac{\theta}{2} \frac{  (\vec{p}\cdot\vec{v})}{|\vec{u}\wedge\vec{v}|^2} .
\end{gather}

Now let us calculate 
\begin{equation}
\mathbf{p}^{\prime\prime} = {\mbox e}^{-\theta\mathcal{I}_2/2} 
\mathbf{p}^\prime {\mbox e}^{\theta\mathcal{I}_2/2} . 
\end{equation}
From eq.\eqref{euler.2} we have 
\begin{equation}
\label{aux.p.prime2}
\begin{split}
\mathbf{p}^{\prime\prime} & = 
 \cos^2\frac{\theta}{2} 
-\cos\frac{\theta}{2}\sin\frac{\theta}{2}[\mathcal{I}_2,\mathbf{p}] 
-\sin^2\frac{\theta}{2}\mathcal{I}_2\mathbf{p}\mathcal{I}_2 \\
& \quad + 
C_u \cos^2\frac{\theta}{2}  
-C_u \cos\frac{\theta}{2}\sin\frac{\theta}{2}[\mathcal{I}_2,\mathbf{u}+\mathbf{u}^\ast] 
-C_u \sin^2\frac{\theta}{2}\mathcal{I}_2(\mathbf{u}+\mathbf{u}^\ast)\mathcal{I}_2 \\
& \quad + 
C_v \cos^2\frac{\theta}{2} 
-C_v \cos\frac{\theta}{2}\sin\frac{\theta}{2}[\mathcal{I}_2,\mathbf{v}+\mathbf{v}^\ast] 
-C_v \sin^2\frac{\theta}{2}\mathcal{I}_2(\mathbf{v}+\mathbf{v}^\ast)\mathcal{I}_2 .
\end{split}
\end{equation}
The calculations are analogous to the previous case, and we 
obtain that 
\begin{equation}
\mathbf{p}^{\prime\prime} = \mathbf{p}^\prime + C_u (\mathbf{u}-\mathbf{u}^\ast) 
+ C_v(\mathbf{v}-\mathbf{v}^\ast) , 
\end{equation}
that is, 
\begin{equation}
\mathbf{p}^{\prime\prime} = \mathbf{p} + 2C_u \mathbf{u} + 
2 C_v \mathbf{v} . 
\end{equation}
If we decompose $\vec{p}$ as 
in eq.\eqref{decomposition.p} and use $|\vec{v}|^2 = |\vec{u}|^2 = 1$, 
we conclude that
\begin{equation}
\mathbf{p}^{\prime\prime} = \mathbf{p}_\perp + 
\mathbf{u}\left[\cos\theta \, p_u + \sin\theta \, p_v\right] + 
\mathbf{v}\left[ \cos\theta \, p_v - \sin\theta \, p_u \right] , 
\end{equation}
that is, we have a rotation by an angle $\theta$ in the plane defined by the
vectors $\vec{u}$ and $\vec{v}$. Note that $\theta$ is 
considered positive when measured from $\vec{v}$ to $\vec{u}$.
\end{proof}

\begin{proof}
To prove Theorem~\ref{thm.5},
let us proceed with $\mathcal{S}$ in eq.\eqref{eq.HS} just as  
we did above with $\mathcal{R}$. 
First we note that we can write 
\begin{equation}
\mathcal{S} = \mathcal{S}_1 - \mathcal{S}_2 ,
\end{equation}
with 
\begin{equation}
\mathcal{S}_1 = \frac{1}{2}[\mathbf{u}-\mathbf{u}^\ast,\mathbf{v}+\mathbf{v}^\ast] , 
\qquad 
\mathcal{S}_2 = \frac{1}{2}[\mathbf{u}+\mathbf{u}^\ast, 
\mathbf{v}-\mathbf{v}^\ast] . 
\end{equation}
Then we have 
\begin{equation}
\label{aux.s.12}
(\mathcal{S}_1)^2 = (\mathcal{S}_2)^2 = |\mathbf{u}|^2 |\mathbf{v}|^2 , 
\end{equation}
and 
\begin{equation}
\label{aux.s.comut}
[\mathcal{S}_1,\mathcal{S}_2] = 4(\vec{v}\cdot\vec{u}) 
(\mathbf{u}\mathbf{v} + \mathbf{u}^\ast\mathbf{v}^\ast) . 
\end{equation}
While eq.\eqref{aux.s.12} shows that we have simple expressions
for ${\mbox e}^{t\mathcal{S}_1}$ and ${\mbox e}^{t\mathcal{S}_2}$, 
eq.\eqref{aux.s.comut} shows that ${\mbox e}^{t(\mathcal{S}_1-\mathcal{S}_2)} 
\neq {\mbox e}^{t\mathcal{S}_1} {\mbox e}^{-t\mathcal{S}_2}$. 
However, the situation can be bypassed if we choose the vectors 
$\vec{u}$ and $\vec{v}$ to be orthogonal,
since from eq.\eqref{aux.s.comut}
\begin{equation}
\vec{v}\cdot\vec{u} = 0 \; \Leftrightarrow \; 
[\mathcal{S}_1,\mathcal{S}_2] = 0 . 
\end{equation}
Let us make this choice and define 
\begin{equation}
\mathcal{H}_1 = \frac{\mathcal{S}_1}{|\vec{u}| |\vec{v}|} , \qquad 
\mathcal{H}_2 = \frac{\mathcal{S}_2}{|\vec{u}| |\vec{v}|} . 
\end{equation}
Then we have 
\begin{equation}
(\mathcal{H}_1)^2 = (\mathcal{H}_2)^2 = 1 ,
\end{equation}
and 
\begin{gather}
{\mbox e}^{\theta\mathcal{H}_1/2} = \cosh\frac{\theta}{2} + 
\mathcal{H}_1 \sinh\frac{\theta}{2} , \\
{\mbox e}^{\theta\mathcal{H}_2/2} = \cosh\frac{\theta}{2} + 
\mathcal{H}_2 \sinh\frac{\theta}{2} .
\end{gather}

The rest of the calculations are completely analogous to the proof of Theorem~\ref{thm.4}, so we will omit the details. We obtain, for 
\begin{equation}
\mathbf{p}^{\prime\prime} = {\mbox e}^{-\theta\mathcal{H}_2/2} 
{\mbox e}^{\theta\mathcal{H}_1/2}\mathbf{p} {\mbox e}^{-\theta\mathcal{H}_1/2} {\mbox e}^{\theta\mathcal{H}_2/2} , 
\end{equation}
that  
\begin{equation}
\mathbf{p}^{\prime\prime} = \mathbf{p}_\perp + 
\mathbf{u}[\cosh\theta p_u + \sinh\theta p_v] + 
\mathbf{v}[\cosh\theta p_v + \sinh\theta p_u] , 
\end{equation}
which we recognize as a hyperbolic rotation by an angle
$\theta$ in the plane of the vectors $\vec{u}$ and
$\vec{v}$.
\end{proof}

\section{Calculation of $(\mathcal{R}_1)^2$ and $(\mathcal{R}_2)^2$}. 
\label{appendix.C}

From eq.\eqref{eq.R1.R2}, the calculation of 
$(\mathcal{R}_1)^2$ and $(\mathcal{R}_2)^2$ involves 
the expression 
\begin{equation}
4 \left(\mathcal{R}_{1|2}\right)^2  =  
\left( 2\mathbf{u}\mathbf{v} + 2\mathbf{u}^\ast \mathbf{v}^\ast 
\pm (\mathbf{u}\mathbf{v}^\ast - \mathbf{v}^\ast\mathbf{u}) 
\pm (\mathbf{u}^\ast\mathbf{v}-\mathbf{v}\mathbf{u}^\ast)\right)^2 , 
\end{equation}
where the upper sign refers to $\mathcal{R}_1$ and the lower
sign refers to $\mathcal{R}_2$. 
When we calculate this product, we obtain a sum 
of terms involving products of four operators. 
The terms where 
any of the operators $\mathbf{u}$, $\mathbf{v}$,  
$\mathbf{v}^\ast$ and $\mathbf{u}$ appears more than once 
can be simplified using the commutation relations 
\eqref{CAR.1.g}, \eqref{CAR.2.g} and \eqref{CAR.3.g}. 
For example: 
\begin{equation}
\mathbf{u}\mathbf{v}\mathbf{v}^\ast \mathbf{u} = 
\mathbf{u}\mathbf{v}((\vec{v}\cdot\vec{u})- 
\mathbf{u}\mathbf{v}) = (\vec{v}\cdot\vec{u}) \mathbf{u}\mathbf{v} . 
\end{equation} 
After doing this simplification with all terms of this kind and after the 
cancellation of some terms, we find that
\begin{equation}
\label{calculation.R.square}
\begin{split}
4 \left(\mathcal{R}_{1|2}\right)^2 & = 
2(\vec{u}\cdot\vec{v})^2 + 4\mathbf{u}\mathbf{v}\mathbf{u}^\ast 
\mathbf{v}^\ast + 4 \mathbf{u}^\ast\mathbf{v}^\ast \mathbf{u}\mathbf{v} + 
\mathbf{u}\mathbf{v}^\ast \mathbf{u}^\ast \mathbf{v} - 
\mathbf{u}\mathbf{v}^\ast \mathbf{v}\mathbf{u}^\ast \\
& \quad -\mathbf{v}^\ast \mathbf{u}\mathbf{u}^\ast\mathbf{v} + 
\mathbf{v}^\ast \mathbf{u}\mathbf{u}\mathbf{v}^\ast + 
\mathbf{u}^\ast\mathbf{v}\mathbf{u}\mathbf{v}^\ast - 
\mathbf{u}^\ast\mathbf{v}\mathbf{v}^\ast \mathbf{u} - 
\mathbf{v}\mathbf{u}^\ast \mathbf{u}\mathbf{v}^\ast + 
\mathbf{v}\mathbf{u}^\ast \mathbf{v}^\ast \mathbf{u} . 
\end{split}
\end{equation}
Note the cancellation of the terms with different signs in the expressions for 
$(\mathcal{R}_1)^2$ and $(\mathcal{R}_2)^2$, so that 
$(\mathcal{R}_1)^2 = (\mathcal{R}_2)^2$. Now let us rearrange all 
these terms to write them in terms of
$\mathbf{u}\mathbf{v}\mathbf{u}^\ast\mathbf{v}^\ast$. 
We have 
\begin{alignat}{1}
& \mathbf{u}^\ast\mathbf{v}^\ast\mathbf{u}\mathbf{v} = 
-|\vec{u}|^2 |\vec{v}|^2 + 
(\vec{v}\cdot\vec{u})(\mathbf{u}^\ast\mathbf{v} - \mathbf{u}
\mathbf{v}^\ast) + |\vec{u}|^2 \mathbf{v}\mathbf{v}^\ast + 
|\vec{v}|^2 \mathbf{u}\mathbf{u}^\ast + \mathbf{u}\mathbf{v}
\mathbf{u}^\ast\mathbf{v}^\ast , \\
& \mathbf{u}\mathbf{v}^\ast \mathbf{u}^\ast \mathbf{v} = 
(\vec{u}\cdot\vec{v})\mathbf{u}\mathbf{v}^\ast - 
|\vec{v}|^2\mathbf{u}\mathbf{u}^\ast  - \mathbf{u}\mathbf{v}
\mathbf{u}^\ast \mathbf{v}^\ast , \\
& \mathbf{u}\mathbf{v}^\ast \mathbf{v}\mathbf{u}^\ast = |\vec{v}|^2 
\mathbf{u}\mathbf{u}^\ast + \mathbf{u}\mathbf{v}\mathbf{u}^\ast \mathbf{v}^\ast , \\
& \mathbf{v}^\ast\mathbf{u}\mathbf{u}^\ast \mathbf{v} = 
(\vec{v}\cdot\vec{u})^2 - (\vec{v}\cdot\vec{u})
(\mathbf{u}\mathbf{v}^\ast + \mathbf{v}\mathbf{u}^\ast) + 
|\vec{v}|^2 \mathbf{u}\mathbf{u}^\ast + 
\mathbf{u}\mathbf{v}\mathbf{u}^\ast \mathbf{v}^\ast , \\
& \mathbf{v}^\ast \mathbf{u}\mathbf{v}\mathbf{u}^\ast = 
(\vec{v}\cdot\vec{u})\mathbf{v}\mathbf{u}^\ast -
|\vec{v}|^2 \mathbf{u}\mathbf{u}^\ast - 
\mathbf{u}\mathbf{v}\mathbf{u}^\ast \mathbf{v}^\ast , \\
& \mathbf{u}^\ast \mathbf{v}\mathbf{u}\mathbf{v}^\ast = 
(\vec{v}\cdot\vec{u})\mathbf{u}\mathbf{v}^\ast -
|\vec{u}|^2 \mathbf{v}\mathbf{v}^\ast - 
\mathbf{u}\mathbf{v}\mathbf{u}^\ast \mathbf{v}^\ast , \\
& \mathbf{u}^\ast \mathbf{v}\mathbf{v}^\ast \mathbf{u} = 
(\vec{v}\cdot\vec{u})^2 - (\vec{v}\cdot\vec{u}) 
(\mathbf{u}\mathbf{v}^\ast + \mathbf{v}\mathbf{u}^\ast) + 
|\vec{u}|^2 \mathbf{v}\mathbf{v}^\ast + 
\mathbf{u}\mathbf{v} \mathbf{u}^\ast \mathbf{v}^\ast , \\
& \mathbf{v}\mathbf{u}^\ast \mathbf{u}\mathbf{v}^\ast = 
|\vec{u}|^2 \mathbf{v}\mathbf{v}^\ast + 
\mathbf{u}\mathbf{v}\mathbf{u}^\ast \mathbf{v}^\ast , \\
& \mathbf{v}\mathbf{u}^\ast \mathbf{v}^\ast \mathbf{u} = 
(\vec{v}\cdot\vec{u}) \mathbf{v}\mathbf{u}^\ast - 
|\vec{u}|^2 \mathbf{v}\mathbf{v}^\ast - 
\mathbf{u}\mathbf{v} \mathbf{u}^\ast \mathbf{v}^\ast . 
\end{alignat}
Using all these expressions in eq.\eqref{calculation.R.square}, 
cancelling terms and using eq.\eqref{CAR.3.g}, we conclude that
\begin{equation}
4 \left(\mathcal{R}_{1|2}\right)^2 = 4(\vec{v}\cdot\vec{u})^2 - 
4|\vec{v}|^2 |\vec{u}|^2 , 
\end{equation}
that is, 
\begin{equation}
\left(\mathcal{R}_{1|2}\right)^2 = - |\vec{u}\wedge\vec{v}|^2 = 
- \left[  
|\vec{v}|^2 |\vec{u}|^2 - 
(\vec{v}\cdot\vec{u})^2\right] < 0 . 
\end{equation}

\end{document}